\documentclass{amsproc}

\usepackage{amsmath}
\usepackage{amssymb}
\usepackage{amsthm}
\usepackage{amscd}
\usepackage{graphicx}

\newcommand{\hako}[1]{\fbox{\parbox[c][6.5pt][b]{4.5pt}{$ #1 $}}}
\newcommand{\hakoc}[1]{\fbox{\parbox[c][6.5pt][b]{6pt}{$ #1 $}}}

\newcommand{\hakoww}[1]{\fbox{\parbox[c][6.5pt][b]{22pt}{$ #1 $}}}

\newcommand{\Z}{{\mathbb Z}}

\newtheorem{theorem}{Theorem}[section]
\newtheorem{lemma}[theorem]{Lemma}

\newtheorem{proposition}[theorem]{Proposition}
\newtheorem{corollary}[theorem]{Corollary}

\theoremstyle{definition}

\newtheorem{example}[theorem]{Example}

\numberwithin{equation}{section}

\makeatletter
\def\eqnarray{%
  \stepcounter{equation}%
  \let\@currentlabel=\theequation
  \global\@eqnswtrue
  \global\@eqcnt\z@
  \tabskip\@centering
  \let\\=\@eqncr
  $$\halign to \displaywidth\bgroup\@eqnsel\hskip\@centering
  $\displaystyle\tabskip\z@{##}$&\global\@eqcnt\@ne
  \hfil$\displaystyle{{}##{}}$\hfil
  &\global\@eqcnt\tw@$\displaystyle\tabskip\z@{##}$\hfil
  \tabskip\@centering&\llap{##}\tabskip\z@\cr}
\makeatother

\newfont{\germ}{eufm10}
\newfont{\slsmall}{cmsl8}

\def\Aff{\mbox{\sl Aff}}
\def\bt{\tilde{b}}
\def\cd{\cdots}
\def\et#1{\tilde{e}_{#1}}
\def\ft#1{\tilde{f}_{#1}}
\def\gehb{\overline{\geh}}

\def\La{\Lambda}
\def\la{\lambda}
\def\Lab{\overline{\La}}
\def\ol#1{\overline{#1}}
\def\ot{\otimes}
\def\Proof{\noindent{\sl Proof.}\quad}
\def\qed{~\rule{1mm}{2.5mm}}

\def\vertex(#1,#2,#3,#4){
\setlength{\unitlength}{1mm}
\begin{picture}(20,20)(0,-5)
\put(3,5){\line(1,0){8}}
\put(7,1){\line(0,1){8}}
\put(0,4.3){$#1$}
\put(12,4.3){$#4$}
\put(6.2,10.5){$\scriptstyle #2$}
\put(6.2,-2.5){$\scriptstyle #3$}
\end{picture}}
\def\btilde{\tilde{b}}
\def\R{{\mathcal R}}
\def\geh{\mathfrak{g}}

\def\wt{\mbox{\sl wt}\,}
\def\cd{\cdots}

\def\veps{\varepsilon}
\def\vphi{\varphi}

\def\Z{{\mathbb Z}}

\newlength{\sukima}
\newlength{\senhaba}
\newlength{\mojitosen}
\newlength{\sentosen}
\newlength{\tatesendepth}
\newcommand{\tatesen}{\vrule height0.4pt width\senhaba depth\tatesendepth 
\hskip-\senhaba}
\newcommand{\yokosen}{\hrule height\senhaba}

\newcommand{\kakomiu}[1]{%
	\vbox{\hbox{\tatesen$#1$\tatesen}
	\vskip-\sentosen \vskip-\senhaba \vskip-\sentosen \vskip-\senhaba%
	\yokosen \vskip-\senhaba \vskip-\mojitosen}}

\newcommand{\kakomid}[1]{%
	\vbox{\hbox{\tatesen$#1$\tatesen}
	\vskip-\sentosen \vskip-\senhaba
	\yokosen
	\vskip-\senhaba \vskip-\sentosen \vskip-\senhaba \vskip-\mojitosen}}

\newcommand{\kakomis}[1]{%
	\vbox{\hbox{\tatesen$#1$\tatesen}
	\yokosen
	\vskip-\senhaba \vskip-\tatesendepth}}

\newcommand{\kakomiud}[1]{%
	\vbox{\hbox{\tatesen$#1$\tatesen}
	\vskip-\sentosen \vskip-\senhaba \vskip-\sentosen \vskip-\senhaba
	\yokosen
	\vskip\sentosen
	\yokosen
	\vskip-\senhaba \vskip-\sentosen \vskip-\senhaba \vskip-\mojitosen}}

\newcommand{\kakomius}[1]{%
	\vbox{\hbox{\tatesen$#1$\tatesen}
	\vskip-\sentosen \vskip-\senhaba \vskip-\sentosen \vskip-\senhaba
	\yokosen
	\vskip\sentosen \vskip\senhaba \vskip\sentosen
	\yokosen
	\vskip-\senhaba \vskip-\tatesendepth}}

\newcommand{\kakomids}[1]{%
	\vbox{\hbox{\tatesen$#1$\tatesen}
	\vskip-\sentosen \vskip-\senhaba
	\yokosen
	\vskip\sentosen
	\yokosen
	\vskip-\senhaba \vskip-\tatesendepth}}

\begin{document}

\title[Scattering rules in soliton cellular automata]
{Scattering rules in soliton cellular automata \\
associated with crystal bases}

\author{G. Hatayama}
\address{Institute of Physics, University of Tokyo, Tokyo 153-8902, Japan}
\email{hata@gokutan.c.u-tokyo.ac.jp}

\author{A. Kuniba}
\address{Institute of Physics, University of Tokyo, Tokyo 153-8902, Japan}
\email{atsuo@gokutan.c.u-tokyo.ac.jp}

\author{M. Okado}
\address{Department of Informatics and Mathematical Science,
    Graduate School of Engineering Science,
Osaka University, Osaka 560-8531, Japan}
\email{okado@sigmath.es.osaka-u.ac.jp}

\author{T. Takagi}
\address{Department of Applied Physics, National Defense Academy,
Kanagawa 239-8686, Japan}
\email{takagi@nda.ac.jp}

\author{Y. Yamada}
\address{Department of Mathematics, Faculty of Science,
Kobe University, Hyogo 657-8501, Japan}
\email{yamaday@math.kobe-u.ac.jp}


\begin{abstract}
Solvable vertex models in a ferromagnetic regime give rise to 
soliton cellular automata at $q=0$.
By means of the crystal base theory, 
we study a class of such automata associated with 
the quantum affine algebra $U_q(\geh_n)$  for 
non exceptional series
$\geh_n = A^{(2)}_{2n-1}, A^{(2)}_{2n},
B^{(1)}_n, C^{(1)}_n, D^{(1)}_n$ and $D^{(2)}_{n+1}$.
They possess a commuting family of time evolutions 
and solitons labeled by crystals of the
smaller algebra $U_q(\geh_{n-1})$.
Two-soliton scattering rule is identified with the 
combinatorial $R$  of $U_q(\geh_{n-1})$-crystals, 
and the multi-soliton scattering is shown 
to factorize into the two-body ones.
\end{abstract}

\maketitle

\section{Introduction}

\subsection{Background}

In early 90's the crystal base theory explained \cite{KMN}
the appearance of  affine Lie algebra characters \cite{DJKMO}
in Baxter's corner transfer matrix method \cite{B}.
The success concerned solvable vertex models in an 
antiferromagnetic regime.
Typically for the six-vertex model corresponding to 
$U_q(A^{(1)}_1)$, the spins in this regime
are alternating at $q=0$.
On the other hand the ferromagnetic regime, where 
the distant spins are all equal, had been left intact
until rather recently.
A renaissance arose with the discovery \cite{HKT1,FOY,HHIKTT,HKT2}
that the ferromagnetic vertex models at $q=0$ turn out to be 
soliton cellular automata \cite{TS,T}.
As a result the ultradiscrete soliton equation \cite{TTMS} 
acquired an origin also as a quantum integrable system at $q=0$.

In \cite{HKT1} a class of soliton cellular automata
were constructed associated with  
$U'_q(\geh_n)$-crystals for non exceptional series 
$\geh_n = A^{(2)}_{2n-1}, A^{(2)}_{2n},
B^{(1)}_n, C^{(1)}_n, D^{(1)}_n$ and $D^{(2)}_{n+1}$.
The solitons therein possess internal degrees of 
freedom labeled by  
crystals of the smaller algebra $U'_q(\geh_{n-1})$.

The purpose of this paper is to study the scattering 
of the solitons in these automata.
The two-body scattering is characterized by the data:
\begin{itemize}
\item[$\bullet$] Exchange of internal degrees of freedom.
\item[$\bullet$] Phase shift.
\end{itemize}
In \cite{FOY} it was shown for $\geh_n = A^{(1)}_n$ that
they fit completely with the combinatorial $R$  
and the associated energy function for 
$U_q(A^{(1)}_{n-1})$-crystals. 
Moreover, the multi-soliton scattering 
factorizes into two-body ones.
See also \cite{TNS}.

The main result of this paper is {\sc Theorem} \ref{th:main},
which extends such results to the 
$\geh_n$ automaton \cite{HKT1} for 
all the non exceptional $\geh_n$ mentioned above.
In these automata  local states take values in the 
crystal of the vector representation.
We will also consider the generalized automata 
corresponding to the analogue of symmetric tensor representations
of $U'_q(\geh_n)$.
Except for the phase part, such a generalization has been studied 
for $A^{(1)}_n$ case in \cite{HHIKTT}.
Our proof of {\sc Theorem} \ref{th:main} 
contains a new ingredient, an operator $T_\natural$,  
which transforms a soliton into another.
To illustrate the idea, let us prove the $A^{(1)}_n$ case 
in the next subsection.

\subsection{\mathversion{bold} $A^{(1)}_n$ case}\label{subsection:Atype}

Let us quickly review necessary notations. For $U'_q(A^{(1)}_n)$-crystal
$B_l$, see Appendix \ref{app:Bl}. 
First 
$\imath_l$ is a map sending  an element 
$(x_1,x_2,\cd,x_n)$ of the $U'_q(A^{(1)}_{n-1})$-crystal $B_l$ to 
	$\hakoww{{n\!+\!1}}^{\,\otimes x_n} \otimes \dots \otimes
	\hako{3}^{\,\otimes x_2} \otimes \hako{2}^{\,\otimes x_1}$
in the $U'_q(A^{(1)}_n)$-crystal $B_1^{\ot l}$. It is known \cite{KN} 
that $\imath_l$
is an isomorphism of $U_q(A_{n-1})$-crystals in the sense that 
$\imath_l\cdot\et{i}=\et{i+1}\cdot\imath_l$ and $\imath_l\cdot\ft{i}
=\ft{i+1}\cdot\imath_l$
for $i=1,2,\cd,n-1$. A one-soliton state of our cellular automaton is described as
\[
\dots\otimes\hako{1}\otimes\dots\otimes\hako{1}\otimes\imath_l(b)\otimes
\hako{1}\otimes\dots\otimes\hako{1}\otimes\dots
\]
for some $b\in B_l$. Such a soliton is said to have length $l$. Under the time
evolution $T_r$ (see Section \ref{subsection:states}), 
the length $l$ soliton moves rightward by $\min(l,r)$ 
steps. One can now define the phase $\gamma$ by counting 
the position $\gamma(t)$ at time $t$ 
{}from the right end as $-\min(l,r)t+\gamma$ under the time 
evolution $T_r$. 
See Section \ref{subsection:scattering} for details.
To the one-soliton $\imath_l(b)$
($b\in B_l$) of phase $\gamma$, we associate an element $z^\gamma b$ of
$\Aff(B_l):=\{z^d b\mid b\in B_l,d\in\Z\}$.

Consider a two-soliton state $z^{\gamma_1}b_1\ot z^{\gamma_2}b_2\in\Aff(B_{l_1})\ot
\Aff(B_{l_2})$ ($l_1>l_2,\gamma_1\gg\gamma_2$) illustrated as
\[
....\imath_{l_1}(b_1)....\imath_{l_2}(b_2)............................................
\]
If $r$ is greater than $l_2$, $\imath_{l_1}(b_1)$ moves faster than 
$\imath_{l_2}(b_2)$ 
under $T_r$. After sufficiently long time, we have
\[
.......................................
......\imath_{l_2}(b'_2)...\imath_{l_1}(b'_1)....
\]
which contains solitons of the same length $l_i$ ($i=1,2$) as the initial state. 
Let 
$\gamma'_i$ be the phase of the soliton $\imath_{l_i}(b'_i)$. 
Now the theorem is
\medskip

\noindent
{\sc Theorem}.\quad 
{\it 
The scattering rule:
\[
z^{\gamma_1}b_1\ot z^{\gamma_2}b_2\longmapsto 
z^{\gamma'_2}b'_2\ot z^{\gamma'_1}b'_1
\]
is given by the combinatorial $R$ :
\[
R:z^{\gamma_1}b_1\ot z^{\gamma_2}b_2\longmapsto 
z^{\gamma_2+\delta}b'_2\ot z^{\gamma_1-\delta}b'_1,
\]
where $\delta=H(b_1\ot b_2)$ is the energy function.
(Therefore $\gamma'_1= \gamma_1 -\delta, \gamma'_2= \gamma_2+\delta$.) }
\medskip

The proof goes as follows. First note that $T_r$ and $R$ both commute with
the actions of $\et{i},\ft{i}$ for $i=1,\cd,n-1$ on $B_{l_1}\ot B_{l_2}$ or
$B_{l_2}\ot B_{l_1}$ ($U_q(A_{n-1})$-invariance). Thus we are allowed to 
reduce the proof to the set of highest weight elements, {\it i.e.}, 
elements killed by $\et{i}$ ($i=1,\cd,n-1$). The set of such elements in
$B_{l_1}\ot B_{l_2}$ is denoted by ${\mathcal H}(B_{l_1}\ot B_{l_2})$. 
If $l_1\ge l_2$,
we have
\begin{eqnarray*}
{\mathcal H}(B_{l_1}\ot B_{l_2})&=&\{u_{l_1}\ot(x_1,x_2,0,\cd,0)
\mid x_1+x_2=l_2\},\\
{\mathcal H}(B_{l_2}\ot B_{l_1})&=&\{u_{l_2}\ot(y_1,y_2,0,\cd,0)\mid 
y_1+y_2=l_1,y_2\le l_2\},
\end{eqnarray*}
where $u_l$ stands for $(l,0,\cd,0)\in B_l$. It is known that 
\begin{eqnarray*}
R:\quad \Aff(B_{l_1})\ot\Aff(B_{l_2})&\longrightarrow&
\Aff(B_{l_2})\ot\Aff(B_{l_1}),\\
z^{\gamma_1}u_{l_1}\ot z^{\gamma_2}(x_1,x_2,0,\cd,0)&\longmapsto&
z^{\gamma_2+H}u_{l_2}\ot z^{\gamma_1-H}(x_1+l_1-l_2,x_2,0,\cd,0),
\end{eqnarray*}
where $H=x_1+l_2$ under suitable normalization of the energy function.
(See (\ref{eq:Hnormalize}).)

For the proof we make use of another operator $T_\natural$ constructed {}from the 
$U'_q(A^{(1)}_n)$-crystal $B_\natural$, which corresponds to the 
degree 2 antisymmetric tensor representation.
$T_\natural$ is to put the highest element ${1\choose2} \in B_\natural$ on the 
left and to carry it through to the right via the 
isomorphism 
$B_\natural \ot B_1 \ot \cd \ot B_1 \simeq 
B_1 \ot \cd \ot B_1 \ot B_\natural$.
(See  \cite{NY} for the algorithm to compute 
$B_\natural \ot B_1 \simeq B_1 \ot B_\natural$.)
It acts on the space of states and 
has the following properties:
\begin{itemize}
\item[$\bullet$] $T_\natural$ commutes with $T_r$.
\item[$\bullet$] On a one-soliton state of length $l$, $T_\natural$ acts as
\[
T_\natural(z^\gamma(x_1,x_2,0,\cd,0))=\left\{
\begin{array}{ll}
z^\gamma(l,0,0,\cd,0))&\mbox{ if }x_2=0\\
z^{\gamma-1}(x_1+1,x_2-1,0,\cd,0))&\mbox{ if }x_2>0.
\end{array}\right.
\]
\end{itemize}
If one cuts the figure given below in the middle, the left and the
right pieces become examples of $x_2 = 0$ and $x_2=1$ 
cases in the above, respectively.

Now suppose our theorem is valid for a particular highest weight element
$z^{\gamma_1}u_{l_1}\ot z^{\gamma_2}(x_1,x_2,0,\cd,0)$ ($x_2>0$), {\it i.e.},
\begin{eqnarray*}
&&T_r^t(z^{\gamma_1}u_{l_1}\ot z^{\gamma_2}(x_1,x_2,0,\cd,0))\\
&&\qquad=z^{\gamma_2+x_1+l_2}u_{l_2}\ot z^{\gamma_1-x_1-l_2}(x_1+l_1-l_2,x_2,0,\cd,0)
\end{eqnarray*}
for sufficiently large $t$.
Apply $T_\natural$ on both sides and use the commutativity 
$T_\natural T_r=T_r T_\natural$. We get
\begin{eqnarray*}
&&T_r^t(z^{\gamma_1}u_{l_1}\ot z^{\gamma_2-1}(x_1+1,x_2-1,0,\cd,0))\\
&&\qquad=z^{\gamma_2+x_1+l_2}u_{l_2}\ot 
z^{\gamma_1-x_1-l_2-1}(x_1+l_1-l_2+1,x_2-1,0,\cd,0).
\end{eqnarray*}
This shows our theorem is also valid for 
$z^{\gamma_1}u_{l_1}\ot z^{\gamma_2-1}(x_1+1,x_2-1,0,\cd,0)$. 
The following figure explains the action of $T_\natural$ 
on a two-soliton state (It transfers the upper line into 
the lower line):
\[
T_\natural(z^{\gamma_1}u_{3}\ot z^{\gamma_2}(1,1,0,\cd,0)) = 
z^{\gamma_1}u_{3}\ot z^{\gamma_2-1}(2,0,0,\cd,0).
\]

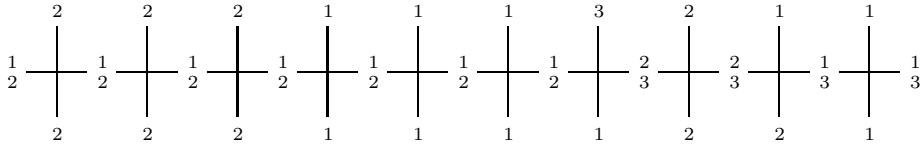
\setlength{\unitlength}{1mm}
\begin{picture}(120,30)(10,0)
\multiput(8,15)(12,0){10}{\line(1,0){8}}
\multiput(12,9)(12,0){10}{\line(0,1){12}}

\put(5,14){${1\atop2}$}
\put(17,14){${1\atop2}$}
\put(29,14){${1\atop2}$}
\put(41,14){${1\atop2}$}
\put(53,14){${1\atop2}$}
\put(65,14){${1\atop2}$}
\put(77,14){${1\atop2}$}
\put(89,14){${2\atop3}$}
\put(101,14){${2\atop3}$}
\put(113,14){${1\atop3}$}
\put(125,14){${1\atop3}$}

\put(11.4,22.3){$\scriptstyle 2$}
\put(23.4,22.3){$\scriptstyle 2$}
\put(35.4,22.3){$\scriptstyle 2$}
\put(47.4,22.3){$\scriptstyle 1$}
\put(59.4,22.3){$\scriptstyle 1$}
\put(71.4,22.3){$\scriptstyle 1$}
\put(83.4,22.3){$\scriptstyle 3$}
\put(95.4,22.3){$\scriptstyle 2$}
\put(107.4,22.3){$\scriptstyle 1$}
\put(119.4,22.3){$\scriptstyle 1$}

\put(11.4,5.9){$\scriptstyle 2$}
\put(23.4,5.9){$\scriptstyle 2$}
\put(35.4,5.9){$\scriptstyle 2$}
\put(47.4,5.9){$\scriptstyle 1$}
\put(59.4,5.9){$\scriptstyle 1$}
\put(71.4,5.9){$\scriptstyle 1$}
\put(83.4,5.9){$\scriptstyle 1$}
\put(95.4,5.9){$\scriptstyle 2$}
\put(107.4,5.9){$\scriptstyle 2$}
\put(119.4,5.9){$\scriptstyle 1$}
\end{picture}

One should notice
that the above procedure is invertible. We thus reduced the theorem to the simplest
case $z^{\gamma_1}u_{l_1}\ot z^{\gamma_2}u_{l_2}$, where the scattering rule 
can be verified directly.

\subsection{Plan of the paper}\label{subsection:plan}

In the rest of the paper we simply write 
$\geh$ to mean $\geh_n$, and let 
$\gehb(= \gehb_n) = C_n, C_n, B_n , C_n, D_n$ and $B_n$  denote 
the classical parts of $\geh = A^{(2)}_{2n-1}, A^{(2)}_{2n},
B^{(1)}_n, C^{(1)}_n, D^{(1)}_n$ and $D^{(2)}_{n+1}$, respectively.
As for the crystals we will use the same symbols $B_l$ and $\Aff(B_l)$ 
either for $\geh$ or $\geh_{n-1}$  as the distinction 
is clear {}from the context.
$B_l$ and $\Aff(B_l)$ will be called 
$U'_q(\geh)$ and $U_q(\geh)$-crystals, respectively.

Our  proof for the listed $\geh$ is basically similar to Section 
\ref{subsection:Atype}.
We start with dividing the cases into three types:
\[
\begin{array}{ll}
\mbox{Type I}& A^{(2)}_{2n-1},B^{(1)}_n,D^{(1)}_n\\
\mbox{Type II}& C^{(1)}_n\\
\mbox{Type III}& A^{(2)}_{2n},D^{(2)}_{n+1}
\end{array}
\]
In each Type the decomposition of $B_l$ or $B_\natural$ into
$U_q(\gehb)$-crystals and the parameterization of the highest weight elements
of $B_{l_1}\ot B_{l_2}$ are the same.
\begin{eqnarray*}
B_l&=&B(l\Lab_1)\quad\mbox{for Type I},\\
&=&B(l\Lab_1)\oplus B((l-2)\Lab_1)\oplus\cd\oplus
\left(B(\Lab_1)\mbox{ or }B(0)\right)\quad\mbox{for Type II},\\
&=&B(l\Lab_1)\oplus B((l-1)\Lab_1)\oplus\cd\oplus B(0)\quad\mbox{for Type III}.\\
B_\natural&=&B(\Lab_2)\oplus B(0)\quad\mbox{for Type I},\\
&=&B(\Lab_2)\quad\mbox{for Type II},\\
&=&B(\Lab_2)\oplus B(\Lab_1)\oplus B(0)\quad\mbox{for Type III}.
\end{eqnarray*}
For Type I and II the proof proceeds just as in the $A^{(1)}_n$ case.
The only difference is that $T_\natural$ is not invertible. But it can
be remedied by virtue of the fact that the weights in the inverse image are
distinct. For Type III the proof is reduced to Type II by using some 
embedding of $U'_q(\geh)$-crystals into $U'_q(C^{(1)}_n)$-crystals.

In Section \ref{section:SCA} we recall the basic feature of 
the $\geh$ automaton. Soliton states and their label in terms of 
tensor products of $\Aff(B_l)$ are explained. 
Our main result is {\sc Theorem} \ref{th:main}.
In Section \ref{sec:proof} we give a proof 
along the idea explained in Section \ref{subsection:Atype}.
In Section \ref{section:examples} typical examples of the 
scattering of solitons are presented.
In Section \ref{sec:discussion} we consider the 
generalization where the automaton states 
are taken {}from $\cd \ot B_{m_j} \ot B_{m_{j+1}} \ot \cd$.
We show that the inhomogeneity of $\{m_j\}$ affects 
the soliton scattering only via a certain phase shift.
Appendix \ref{section:CCR} is a summary of basic facts 
on crystals and combinatorial $R$ .
Appendix \ref{app:Bl} and \ref{section:B-natural} 
contain some information on the crystals 
$B_l$ and $B_\natural$, respectively.

We close the section with  a remark.
Solitons  in our automata with an equal velocity can form a 
bound state (breather) in general.
It propagates stably if isolated, and  
exhibits curious reactions with other solitons or breathers.
In this paper we shall only study the 
scattering of  solitons with  
distinct velocities, and leave the study of bound states
as a future problem.

\section{\mathversion{bold}Soliton cellular automata}
\label{section:SCA}

\subsection{\mathversion{bold}States and time evolutions}
\label{subsection:states}

For $\geh = A^{(2)}_{2n-1} (n \ge 3), A^{(2)}_{2n} (n \ge 2), 
B^{(1)}_n (n \ge 3),
D^{(1)}_n (n \ge 4), D^{(2)}_{n+1} (n \ge 2)$, let $B_l$ be the 
$U'_q(\geh)$-crystal in \cite{KKM}.
For $\geh = C^{(1)}_{n} (n \ge 2)$ 
our $B_l$ is the one in \cite{HKKOT,HKOT1}.
We have listed the parameterization of $B_l$ as a set in Appendix \ref{app:Bl}.
For any $\geh$ we write
\begin{equation}\label{eq:udef}
u_l = (l,0,\ldots,0) \in B_l
\end{equation}
in the notation therein. 

Consider the crystal $B_1^{\ot N}$ for sufficiently large $N$.
The elements of $B_1^{\ot N}$ we have in mind are of the 
following form:
\begin{equation*}
\cd \ot \hako{1} \ot \cd \ot \hako{1} \ot 
b_1 \ot \cd \ot b_k \ot \hako{1} \ot \cd \ot \hako{1} \ot \cd,
\end{equation*}
where $b_1, \ldots, b_k \in B_1$ and $\hako{1} = u_1$.
Namely, relatively few elements are non $\hako{1}$,
and almost all are $\hako{1}$.
See Appendix \ref{app:Bl} for the definition of the symbol $\hako{a}$.
In the assertions below, we embed, if necessary, 
$B_1^{\ot N}$ into $B_1^{\ot N'} (N < N')$ by
\begin{eqnarray*}
B_1^{\ot N}&\hookrightarrow& B_1^{\ot N'}\\
c_1\ot\cd\ot c_N&\mapsto& c_1\ot\cd\ot c_N\ot
\underbrace{\hako{1}\ot\cd\ot \hako{1}}_{N'-N}.
\end{eqnarray*}

Let $B_l \ot B_1 \simeq B_1 \ot B_l$ be an isomorphism 
explained in Appendix \ref{section:CCR}.
\begin{lemma} \label{lem:1}
By iterating  $B_l\ot B_1\simeq B_1\ot B_l$, we consider a map
\begin{eqnarray*}
B_l\ot B_1\ot\cd\ot B_1&\mathop{\longrightarrow}^\sim&
B_1\ot\cd\ot B_1\ot B_l\\
u_l\ot b_1\ot\cd\ot b_N&\mapsto&\btilde_1\ot\cd\ot\btilde_N\ot\btilde,
\end{eqnarray*}
where we assume $b_j = \hako{1}$ for $N' \le j \le N$.
Then there exists an integer $N_0$ such that $\btilde=u_l$ for $N-N' \ge N_0$.
\end{lemma}
Taking sufficiently large $N-N'$ such that the above lemma holds, we define
a map $T_l:B_1^{\ot N}\longrightarrow B_1^{\ot N}$ by 
$b_1\ot\cd\ot b_N\mapsto\btilde_1\ot\cd\ot\btilde_N$.

\begin{lemma} \label{lem:2}
For a fixed element of $B_1^{\ot N}$ as in {\sc Lemma} \ref{lem:1},
there exists an integer $l_0$ such that
$T_l=T_{l_0}$ for any $l\ge l_0$.
\end{lemma}
Both lemmas can be verified directly by using the result in \cite{HKOT1,HKOT2}.
In particular $m=1$ case of the fact
\begin{equation}\label{eq:uu}
u_l \ot u_m \simeq u_m \ot u_l \; \text{ under }\;
B_l \ot B_m \simeq B_m \ot B_l \; \text{for any } l, m,
\end{equation}
is relevant here.

An element of $B_1^{\ot N}$ having the property described in 
the beginning of this subsection will be called a {\it state}. 
{\sc Lemma} \ref{lem:1} and \ref{lem:2} enable
us to define an operator $T=\lim_{l\rightarrow\infty}T_l$ on the space of 
states.
Application of $T$ induces a transition of state. Thus it can be 
regarded as a certain dynamical system, in which $T$ plays
the role of `time evolution'. 
We call it $\geh$ automaton.
By the same reason, $T_l$ may also be viewed as another time evolution.
(In this paper, time evolution means the one by $T$ unless otherwise stated.)

Under the condition stated after {\sc Lemma} \ref{lem:1},
the sequence  $b_1\ot\cd\ot b_N$ is determined uniquely {}from
$\btilde_1\ot\cd\ot\btilde_N$ and $\btilde=u_l$, since the
isomorphism $B_l \ot B_1 \simeq B_1 \ot B_l$ is bijective.
Hence, the time evolutions $T_l$ ($l \geq 1$) are invertible.

\subsection{\mathversion{bold}Conservation laws}
\label{subsection:conservation}

Fix sufficiently large $N$ and consider a composition of the 
combinatorial $R$'s
\[
\R_l=R_{N\,N+1}\cd R_{23}R_{12}:
\Aff(B_l)\ot\Aff(B_1)^{\ot N}\longrightarrow\Aff(B_1)^{\ot N}\ot\Aff(B_l).
\]
Here $R_{i\,i+1}$ signifies that the $R$  acts on the $i$-th and
$(i+1)$-th components of the tensor product.
Applying $\R_l$ to an element $u_l \ot p$ ($p=b_1\ot\cd\ot b_N$), we
have 
\begin{eqnarray*}
\R_l(u_l\ot p)
&=&
z^{H_1}\btilde_1\ot 
z^{H_2}\btilde_2\ot\cd\ot
z^{H_N}\btilde_N\ot z^{\bar{E}_l(p)}u_l,\\
\bar{E}_l(p)&=&-\sum_{j=1}^N H_j, \quad
H_j=H(b^{(j-1)}\ot b_j),
\end{eqnarray*}
where $b^{(0)} = u_l$ and $b^{(j)}$ ($1\le j<N$) is defined by
\begin{eqnarray*}
B_l\ot\underbrace{B_1\ot\cd\ot B_1}_j&\simeq&\underbrace{B_1\ot\cd\ot B_1}_j\ot B_l\\
u_l\ot b_1\ot\cd\ot b_j&\mapsto&
\btilde_1\ot\cd\ot\btilde_j\ot b^{(j)}.
\end{eqnarray*}


\setlength{\unitlength}{1mm}
\begin{picture}(80,28)(-27,0)

\put(-1,14){$\scriptstyle z^0u_l$}

\put(9,5){$\scriptstyle z^{H_1}{\tilde b}_1$}

\put(19,5){$\scriptstyle z^{H_2}{\tilde b}_2$}

\put(38,5){$\scriptstyle z^{H_N}{\tilde b}_N$}

\put(9,22.7){$\scriptstyle z^{0}{b}_1$}
\put(19,22.7){$\scriptstyle z^{0}{b}_2$}
\put(38,22.7){$\scriptstyle z^{0}{b}_N$}

\put(48,14){$\scriptstyle z^{\bar{E}_l(p)}u_l$}

\put(7,15){\line(1,0){20}}
\put(12,9){\line(0,1){12}}

\put(29,14.5){$\cdots$}

\put(36,15){\line(1,0){10}}

\put(22,9){\line(0,1){12}}
\put(41,9){\line(0,1){12}}

\end{picture}

Since $H(u_l\ot u_1)= 2\varsigma \neq 0$ (see (\ref{eq:Hnormalize})),
the quantity  $\bar{E}_l(p)$ grows linearly with $N$.
To avoid it we introduce a regularized version:
\begin{equation*}
E_l(p) = -\sum_{j=1}^N(H_j - H(u_l\ot u_1)) = \bar{E}_l(p) + NH(u_l \ot u_1),
\end{equation*}
which becomes independent of $N$ for   $N$  sufficiently large.

The Yang-Baxter equation in {\sc Proposition} \ref{prop:YBeq}
implies that as maps {}from $\Aff(B_l)\ot\Aff(B_{l'})\ot\Aff(B_1)^{\ot N}$
to $\Aff(B_1)^{\ot N}\ot\Aff(B_{l'})\ot\Aff(B_l)$, the equality
\begin{equation}
\R_{l'}\R_lR_{12}=R_{N+1\,N+2}\R_l\R_{l'} \label{eq:commuting}
\end{equation}
is valid.
Here $\R_l$ (resp. $\R_{l'}$) acts on $\Aff(B_l)\ot\Aff(B_1)^{\ot N}$
(resp. $\Aff(B_{l'})\ot\Aff(B_1)^{\ot N}$), and acts on the other component
as identity.

\begin{proposition}\label{pr:T_l}
For an element $p\in B_1^{\ot N}$, we have
\begin{itemize}
\item[(1)] $T_lT_{l'}(p)=T_{l'}T_l(p)$.
\item[(2)] $E_l(T_{l'}(p))=E_l(p)$. In particular, $E_l(T(p))=E_l(p)$.
\end{itemize}
\end{proposition}

\Proof
We follow the argument in Section 3.2 of \cite{FOY}.
Consider the element 
$z^0 u_l \ot z^0 u_{l'} \ot p 
\in \Aff(B_l) \ot \Aff(B_{l'}) \ot \Aff(B_1)^{\ot N}$.
Apply the both sides of (\ref{eq:commuting}) to it.
Graphically one gets ($\delta=H(u_l \ot u_{l'})$)


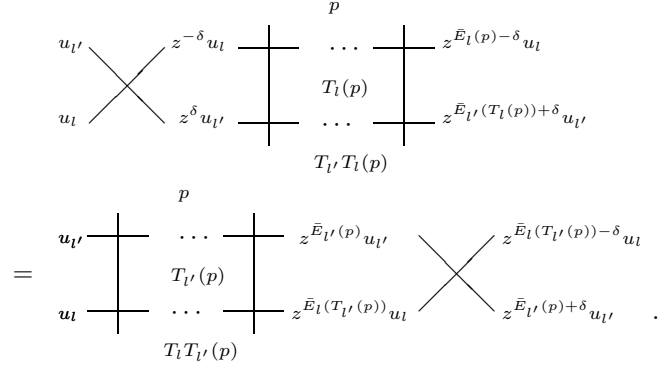
\setlength{\unitlength}{1mm}
\begin{picture}(80,55)(-20,0)

\put(0.2,35){\line(1,1){10}}
\put(0.2,45){\line(1,-1){10}}

\put(-4,35){$\scriptstyle u_{l}$}
\put(-4,45){$\scriptstyle u_{l'}$}

\put(12,35){$\scriptstyle z^{\delta}u_{l'}$}
\put(11,45){$\scriptstyle z^{-\delta}u_{l}$}

\put(20,35){\line(1,0){8}}
\put(20,45){\line(1,0){8}}

\put(31,34){$\cdots$} 
\put(32,44){$\cdots$}

\put(38,35){\line(1,0){8}}
\put(38,45){\line(1,0){8}}

\put(24,32){\line(0,1){16}}
\put(42,32){\line(0,1){16}}

\put(-4,19){$\scriptstyle u_{l'}$}
\put(-4,9){$\scriptstyle u_{l}$}

\put(32,50){$\scriptstyle p$}
\put(31,39){$\scriptstyle T_{l}(p)$}
\put(30,29){$\scriptstyle T_{l'}T_{l}(p)$}

\put(47,45){$\scriptstyle z^{\bar{E}_{l}(p)-\delta}u_{l}$}
\put(47,35){$\scriptstyle z^{\bar{E}_{l'}(T_{l}(p))+\delta}u_{l'}$}


\put(-10,14){$=$}

\put(0,10){\line(1,0){8}}
\put(0,20){\line(1,0){8}}

\put(11,9){$\cdots$} 
\put(12,19){$\cdots$}

\put(18,10){\line(1,0){8}}
\put(18,20){\line(1,0){8}}

\put(4,7){\line(0,1){16}}
\put(22,7){\line(0,1){16}}

\put(-4,19){$\scriptstyle u_{l'}$}
\put(-4,9){$\scriptstyle u_{l}$}

\put(12,25){$\scriptstyle p$}
\put(11,14){$\scriptstyle T_{l'}(p)$}
\put(10,4){$\scriptstyle T_lT_{l'}(p)$}

\put(28,19){$\scriptstyle z^{\bar{E}_{l'}(p)}u_{l'}$}
\put(27,9){$\scriptstyle z^{\bar{E}_{l}(T_{l'}(p))}u_{l}$}

\put(44,10){\line(1,1){10}}
\put(44,20){\line(1,-1){10}}

\put(55,9){$\scriptstyle z^{\bar{E}_{l'}(p)+\delta}u_{l'}$}
\put(55,19){$\scriptstyle z^{\bar{E}_{l}(T_{l'}(p))-\delta}u_{l}$}

\put(75,9){.}

\end{picture}

\noindent
Since $E_l(p)$ differs $\bar{E}_l(p)$ just by a constant,
the assertions (1) and (2) follow directly {}from this diagram.
\qed

In addition to $E_l$, there is another conserved quantity 
resulting {}from the $U_q(\gehb_{n-1})$-invariance explained 
in the beginning of Section \ref{sec:proof}.
However we do not detail it in this paper.
For $A^{(1)}_n$ case, see Section 4.2 of \cite{FOY}.

\subsection{\mathversion{bold}Solitons}
\label{subsection:solitons}

In order to describe solitons in the $\geh$ automaton,
we introduce an injection
$$\imath_l: B_l \text{ for } U'_q(\geh_{n-1}) 
\rightarrow \left(B_1 \text{ for } U'_q(\geh) \right)^{\ot l}$$
for each $l \in {\mathbb Z}_{\ge 1}$ as follows.
%

\noindent
$\geh=A^{(2)}_{2n-1}:$
For $b=(x_1,\dots,x_{n-1},\ol{x}_{n-1},\dots,\ol{x}_1)$
in $U'_q(\geh_{n-1})$-crystal $B_{l}$,
\begin{equation*}
\imath_l (b)=
	\hako{\ol{2}}^{\,\otimes \ol{x}_1} \otimes \dots \otimes
	\hako{\ol{n}}^{\,\otimes \ol{x}_{n\! -\! 1}} \otimes
	\hako{n}^{\,\otimes x_{n\! -\! 1}} \otimes \dots \otimes
	\hako{2}^{\,\otimes x_1}.
\end{equation*}
\noindent
$\geh=A^{(2)}_{2n}:$
For $b=(x_1,\dots,x_{n-1},\ol{x}_{n-1},\dots,\ol{x}_1)$
in $U'_q(\geh_{n-1})$-crystal $B_{l}$,
we define $s(b)=\sum_{i=1}^{n-1}(x_i+\ol{x}_i),\,
s'(b)=[ (l-s(b))/2 ]$.
\vspace*{10pt}\\
\noindent
If $l-s(b)$ is odd,
\[
\imath_l (b)=
	\phi \otimes \hako{\ol{1}}^{\,\otimes s'(b)} \otimes
	\hako{\ol{2}}^{\,\otimes \ol{x}_1} \otimes \dots \otimes
	\hako{\ol{n}}^{\,\otimes \ol{x}_{n\! -\! 1}} \otimes
	\hako{n}^{\,\otimes x_{n\! -\! 1}} \otimes \dots \otimes
	\hako{2}^{\,\otimes x_1} \otimes
	\hako{1}^{\,\otimes s'(b)},
\]
otherwise
\[
\imath_l (b)=
	\hako{\ol{1}}^{\,\otimes s'(b)} \otimes
	\hako{\ol{2}}^{\,\otimes \ol{x}_1} \otimes \dots \otimes
	\hako{\ol{n}}^{\,\otimes \ol{x}_{n\! -\! 1}} \otimes
	\hako{n}^{\,\otimes x_{n\! -\! 1}} \otimes \dots \otimes
	\hako{2}^{\,\otimes x_1} \otimes
	\hako{1}^{\,\otimes s'(b)}.
\]

\noindent
$\geh=B^{(1)}_{n}:$
For $b=(x_1,\dots,x_{n-1},x_0,\ol{x}_{n-1},\dots,\ol{x}_1)$
in $U'_q(\geh_{n-1})$-crystal $B_{l}$,
\[
\imath_l (b)=
	\hako{\ol{2}}^{\,\otimes \ol{x}_1} \otimes \dots \otimes
	\hako{\ol{n}}^{\,\otimes \ol{x}_{n\! -\! 1}} \otimes
	\hako{0}^{\,\otimes x_0} \otimes
	\hako{n}^{\,\otimes x_{n\! -\! 1}} \otimes \dots \otimes
	\hako{2}^{\,\otimes x_1}.
\]
\noindent
$\geh=C^{(1)}_{n}:$
For $b=(x_1,\dots,x_{n-1},\ol{x}_{n-1},\dots,\ol{x}_1)$
in $U'_q(\geh_{n-1})$-crystal $B_{l}$,
we define $s(b)=\sum_{i=1}^{n-1}(x_i+\ol{x}_i),\,
s'(b)=(l-s(b))/2$.
\[
\imath_l (b)=
	\hako{\ol{1}}^{\,\otimes s'(b)} \otimes
	\hako{\ol{2}}^{\,\otimes \ol{x}_1} \otimes \dots \otimes
	\hako{\ol{n}}^{\,\otimes \ol{x}_{n\! -\! 1}} \otimes
	\hako{n}^{\,\otimes x_{n\! -\! 1}} \otimes \dots \otimes
	\hako{2}^{\,\otimes x_1} \otimes
	\hako{1}^{\,\otimes s'(b)}.
\]
\noindent
$\geh=D^{(1)}_{n}:$
For $b=(x_1,\dots,x_{n-1},\ol{x}_{n-1},\dots,\ol{x}_1)$
in $U'_q(\geh_{n-1})$-crystal $B_{l}$,
\[
\imath_l (b)=
	\hako{\ol{2}}^{\,\otimes \ol{x}_1} \otimes \dots \otimes
	\hako{\ol{n}}^{\,\otimes \ol{x}_{n\! -\! 1}} \otimes
	\hako{n}^{\,\otimes x_{n\! -\! 1}} \otimes \dots \otimes
	\hako{2}^{\,\otimes x_1}.
\]
\noindent
$\geh=D^{(2)}_{n+1}:$
For $b=(x_1,\dots,x_{n-1},x_0,\ol{x}_{n-1},\dots,\ol{x}_1)$
in $U'_q(\geh_{n-1})$-crystal $B_{l}$,
we define $s(b)=x_0+\sum_{i=1}^{n-1}(x_i+\ol{x}_i),\,
s'(b)=[ (l-s(b))/2 ]$.
\vspace*{10pt}\\
\noindent
If $l-s(b)$ is odd,
\[
\imath_l (b)=
	\phi \otimes \hako{\ol{1}}^{\,\otimes s'(b)} \otimes
	\hako{\ol{2}}^{\,\otimes \ol{x}_1} \otimes \dots \otimes
	\hako{\ol{n}}^{\,\otimes \ol{x}_{n\! -\! 1}} \otimes
	\hako{0}^{\,\otimes x_0} \otimes
	\hako{n}^{\,\otimes x_{n\! -\! 1}} \otimes \dots \otimes
	\hako{2}^{\,\otimes x_1} \otimes
	\hako{1}^{\,\otimes s'(b)},
\]
otherwise
\[
\imath_l (b)=
	\hako{\ol{1}}^{\,\otimes s'(b)} \otimes
	\hako{\ol{2}}^{\,\otimes \ol{x}_1} \otimes \dots \otimes
	\hako{\ol{n}}^{\,\otimes \ol{x}_{n\! -\! 1}} \otimes
	\hako{0}^{\,\otimes x_0} \otimes
	\hako{n}^{\,\otimes x_{n\! -\! 1}} \otimes \dots \otimes
	\hako{2}^{\,\otimes x_1} \otimes
	\hako{1}^{\,\otimes s'(b)}.
\]

The $\imath_l$ here was originally introduced in \cite{HKT1}, where 
the following was shown.
\begin{lemma}\label{lem:i}
For $i \in \{1,2,\dots,n-1\}$, the following diagram is commutative:
\[
\begin{CD}
U'_q(\geh_{n-1})\text{-crystal }B_l @>{\imath_l}>>
U'_q(\geh)\text{-crystal }B_1^{\otimes l}\\
@V{\tilde{e}_i}VV @VV{\tilde{e}_{i+1}}V \\
U'_q(\geh_{n-1})\text{-crystal }B_l \sqcup\{0\} @>>{\imath_l}>
U'_q(\geh)\text{-crystal }B_1^{\otimes l}\sqcup\{0\},\\
\end{CD}
\]
where we set $\imath_l(0) = 0$.
The same relation holds between $\ft{i}$ and $\ft{i+1}$.
\end{lemma}

The Kashiwara operators of
$U'_q(\geh_{n-1})$ and $U'_q(\geh)$-crystals
should not be confused although we use the same notation.
The proof is due to the explicit
rules for $\tilde{e}_i$ \cite{KKM, HKKOT} and
the embedding of $U'_q(\geh_{n-1})\text{-crystal }B_l$ into
$U'_q(\geh_{n-1})\text{-crystal }B_1^{\ot l}$ as
$U_q(\gehb_{n-1})$-crystals \cite{KN}.

A state of the form 
\begin{equation}\label{eq:soliton-state}
...[l_1]........[l_2]..... \cdots .....[l_m]....... 
\end{equation}
is called an $m$-soliton state of
length $l_1,l_2,\ldots,l_m$.
Here $..[l]..$ denotes a local configuration such as
\[
\cd\ot 1 \ot 1 \ot \imath_l(b) \ot 1 \ot 1 \ot\cd
\quad \mbox{ for some } b \in B_l
\]
surrounded by sufficiently many $1$'s.
In view of (\ref{eq:shift}) below, we assume that there are 
at least $l_i$ 1's (i.e., $1^{\ot l_i}$) between $[l_i]$ and $[l_{i+1}]$.
(Here and in what follows, we often abbreviate 
$\hako{a} \in B_1$ to $a$.)
In the sequel we set 
\begin{equation}\label{eq:varsigma}
\varsigma = \begin{cases}
2 & \mbox{ if } \geh = A^{(2)}_{2n}, D^{(2)}_{n+1},\\
1 & \mbox{ otherwise}. 
\end{cases}
\end{equation}
By combining the results in \cite{HKOT1,HKOT2,HKT1} it is not difficult 
to show
\begin{lemma} \label{lem:single-soliton}
Let $p$ be a one-soliton state of length $l$, then
\begin{itemize}
\item[(1)]
The $k$-th conserved quantity of $p$ is given by $E_k(p)=\varsigma \min(k,l)$.
\item[(2)]
The state $T_k(p)$ is obtained by the rightward shift by
$\min(k,l)$ lattice steps.
\end{itemize}
\end{lemma}

For example (2) is shown by establishing 
\begin{equation}\label{eq:shift}
u_k \ot \imath_l(b)\ot 1^{\ot \min(k,l)} \simeq 
1^{\ot \min(k,l)}\ot \imath_l(b) \ot u_k,
\end{equation}
for any $b \in U'_q(\geh_{n-1})$-crystal $B_l$ 
under the isomorphism $B_k \ot B_1^{\ot l+\min(k,l)} \simeq 
B_1^{\ot l+\min(k,l)} \ot B_k$.

For any state $p$, define the numbers $N_l=N_l(p)$ ($l=1,2,\ldots$) by 
\begin{eqnarray*}
&&E_l=\varsigma\sum_{k \geq 1} \min(k,l) N_k, \quad E_0=0, \\
&&N_l=(-E_{l-1}+2 E_{l}-E_{l+1})/\varsigma.
\end{eqnarray*}
By {\sc Lemma} \ref{lem:single-soliton}, we have
\begin{proposition}
For the $m$-soliton state (\ref{eq:soliton-state}),
$N_l$ is the number of solitons of length $l$, i.e., 
$N_l=\sharp \{ j \ \vert \  l_j=l \}$.
\end{proposition}
This proposition implies the stability of solitons,
since the numbers $E_l(p)$, hence $N_l(p)$, are conserved.

\subsection{Scattering of solitons}\label{subsection:scattering}

We introduce a labeling of solitons of length $l$ using
$\Aff(B_l)$ for the lower rank algebra $U'_q(\geh_{n-1})$.
Suppose there is a soliton of length $l$ at time $t$
\begin{equation}\label{eq:sol}
\cd \ot 1\ot \imath_l(b) \ot 1 \ot \cd
= \cd \ot 1 \ot \nu_1 \ot \nu_2 \ot \cdots \ot \nu_l \ot 1 \ot \cd ,
\end{equation}
where $\nu_1,\ldots, \nu_l \in B_1$ and 
$b$ is an element of $U'_q(\geh_{n-1})$-crystal $B_l$.
Say it is at position $\gamma(t)$, if 
$\nu_1$ is in the $\gamma(t)$-th tensor component of $B_1^{\ot N}$
counting {}from the right end.
{}From {\sc Lemma} \ref{lem:single-soliton} (2), the position $\gamma(t)$ 
under the time evolution $T_k$ is
given by $\gamma(t)=-\min(k,l)t+\gamma$ 
unless it interacts with other solitons. 
Here $\min(k,l)$ is the velocity and 
$\gamma$ is the phase.
To such a soliton (\ref{eq:sol}) we associate the following 
element in  $\Aff(B_l)$.
\begin{alignat}{2}
&z^{\gamma}b &\quad \mbox{if } \geh &\neq A^{(2)}_{2n},D^{(2)}_{n+1},\label{eq:coef1}\\
&z^{2\gamma}b &\quad \mbox{if } \geh &= A^{(2)}_{2n},D^{(2)}_{n+1},
\; \nu_1 \neq \phi,\label{eq:coef2-1}\\
&z^{2\gamma-1}b &\quad \mbox{if } \geh &= A^{(2)}_{2n},D^{(2)}_{n+1},
\; \nu_1 = \phi,\label{eq:coef2-2}
\end{alignat}
where the symbol $\phi \in B_1$ has been 
introduced in Appendix \ref{app:Bl}.

Now consider an $m$-soliton state  illustrated as below.
\[
...[l_1].....[l_2]...............[l_m]................
\]
As in (\ref{eq:soliton-state}) 
we assume that solitons are enough separated so that 
there are at least $l_i$ 1's 
between $[l_i]$ and $[l_{i+1}]$.
To such a state  we associate 
\begin{equation}\label{eq:associate}
z^{\gamma_1}b_1\ot z^{\gamma_2}b_2\ot\cd\ot z^{\gamma_m}b_m
\in \Aff(B_{l_1}) \ot \Aff(B_{l_2}) \ot \cd \ot \Aff(B_{l_m}),
\end{equation}
where each element $z^{\gamma_j}b_j \in \Aff(B_{l_j})$ 
is specified as in the previous paragraph.
Suppose $l_1>l_2>\cd>l_m$. Since longer solitons move faster,
we can expect that the state turns out to be 
\[
.......[l_m]...........[l_2].....[l_1]................
\]
after sufficiently many time evolutions. 
(This claim  is due to {\sc Corollary} \ref{co:main} below.)
We represent such a scattering process as
\begin{eqnarray}
&&z^{\gamma_1}b_1\ot z^{\gamma_2}b_2\ot\cd\ot z^{\gamma_m}b_m \nonumber\\
&\mapsto&z^{\gamma'_m}b'_m\ot\cd\ot z^{\gamma'_2}b'_2\ot z^{\gamma'_1}b'_1.
\label{eq:multi}
\end{eqnarray}

The main result of this paper is the following.
\begin{theorem}\label{th:main}
The two-body scattering of solitons $\imath_l(b_1)$ and
$\imath_k(b_2)$ ($l > k$) under the time evolution $T_r$  ($r>k$)
is described by the combinatorial $R$ :
\begin{eqnarray*}
\Aff(B_{l})\ot \Aff(B_{k})&\simeq&\Aff(B_{k})\ot \Aff(B_{l})\\
z^{m_1}b_1 \ot z^{m_2}b_2 &\mapsto& 
z^{m_2+\delta}b'_2 \ot z^{m_1-\delta}b'_1,
\end{eqnarray*}
with $\delta=H(b_1 \ot b_2)$.
\end{theorem}

The statement on the phase shift is valid 
only for the energy function normalized as (\ref{eq:Hnormalize}).
Such values will be given in (2) of {\sc Propositions} 
\ref{prop:hweI}, \ref{prop:hweII} and {\sc Theorem} \ref{th:RforAD}.

{}From {\sc Theorem} \ref{th:main} 
one can easily establish the following with the help of 
{\sc Proposition} \ref{pr:T_l} and 
{\sc Lemma} \ref{lem:single-soliton} by the same argument 
as in \cite{FOY}.
\begin{corollary}\label{co:main}
Scattering of solitons is factorized into
two-body scattering. Namely, (\ref{eq:multi}) is valid 
under the isomorphism of crystals
\[
\Aff(B_{l_1}) \ot  \cd \ot \Aff(B_{l_m}) \simeq 
\Aff(B_{l_m}) \ot  \cd \ot \Aff(B_{l_1}).
\]
\end{corollary}

The map  decomposes into the combinatorial $R$'s corresponding to 
pairwise transpositions of the components.
It is independent of the order of the transpositions due to the Yang-Baxter
equation ({\sc Proposition} \ref{prop:YBeq}).

Although we have assumed  $l_1 > \cdots > l_m$ in {\sc Corollary} \ref{co:main}, 
some part of the results can be
generalized to the situation 
where there exist solitons with the same length. 
See the second remark after Theorem 4.6 in \cite{FOY}.
However a thorough study of such cases will necessarily involve 
the bound states as mentioned in the end of Section \ref{subsection:plan}.

\section{\mathversion{bold} Proof of {\bf Theorem} \ref{th:main}}
\label{sec:proof}

\subsection{Preliminaries}\label{subsection:pre}

First we recall an important consequence derived {}from representation theory.
Note that $U'_q(\geh)$ contains $U_q(\gehb_{n-1})$
as a subalgebra.
This fact can be translated into the language of crystals,
and combined with {\sc Lemma} \ref{lem:i}, it leads to 
the property that on an $m$-soliton state
\[
T_k\mbox{ commutes with }\et{i},\ft{i}\;(i=1,2,\cd,n-1)\mbox{ on }
\Aff(B_{l_1})\ot\cd\ot\Aff(B_{l_m}).
\]
Here we are identifying the $m$-soliton state with an
element {}from $\Aff(B_{l_1}) \ot \cd \ot \Aff(B_{l_m})$
as explained in Section \ref{subsection:scattering}.
The action of $T_k$ on $\Aff(B_{l_1}) \ot \cd \ot \Aff(B_{l_m})$ is 
defined through this identification. 
The actions of $\et{i},\ft{i}$ ($i=1,2,\cd,n-1$)
are calculated according to the rule explained 
in Appendix \ref{section:CCR}. 
By the actions, the power of $z$ in an element
of $\Aff(B_{l_j})$ is unaffected. We call this property 
$U_q(\gehb_{n-1})$-invariance.
This reduces the proof of {\sc Theorem} \ref{th:main} to 
the subset of $\Aff(B_l) \ot \Aff(B_k)$ corresponding to the 
highest weight elements with respect to $U_q(\gehb_{n-1})$.

So we give the data on the highest weight elements of $B_l\ot B_k$ and their 
image under the combinatorial $R$.
For elements in $B_l$ of a particular form, we employ the following notation:
\begin{equation}\label{eq:notation}
(x_1,x_2,\ol{x}_1)=(x_1,x_2,0,\cd,0,\ol{x}_1),\qquad (x_1)=(x_1,0,\cd,0).
\end{equation}

For Type I case, we have 
\begin{proposition}[\cite{HKOT1,HKOT2}] \label{prop:hweI}
Assume $l\ge k$.
\begin{itemize}
\item[(1)] The set of highest weight elements of $B_l\ot B_k$ is given by
           $\{(l)\ot(d,c,b)\mid b+c+d=k\}$.
\item[(2)] The combinatorial $R$ acts on the set of highest weight elements as
\begin{equation} \label{eq:R-hwe}
R: z^m(l)\ot z^{m'}(d,c,b)\mapsto z^{m'+H}(k)\ot z^{m-H}(l-b-c,c,b),
\end{equation}
where $H=2k-2b-c$.
\end{itemize}
\end{proposition}

For Type II case, we have

\begin{proposition}[\cite{HKOT1}] \label{prop:hweII}
Assume $l\ge k$.
\begin{itemize}
\item[(1)] The set of highest weight elements of $B_l\ot B_k$ is given by
           $\{(f)\ot(d,c,b)\mid f\le l,f\equiv l\,(\mbox{mod }2),
            b+c+d\le k,b+c+d\equiv k\,(\mbox{mod }2),f\ge b+c\}$.
\end{itemize}
For a highest weight element $(f)\ot(d,c,b)$, define nonnegative integers
$a,e$ by $l=2e+f,k=2a+b+c+d$.
\begin{itemize}
\item[(2)] The combinatorial $R$ acts on the set of highest weight elements as
           follows:
\begin{itemize}
\item[(i)] If $a\ge e$,
\[
R: z^m(f)\ot z^{m'}(d,c,b)\mapsto z^{m'+H}(k-2e)\ot z^{m-H}(d+l-k-y,c,b-y),
\]
where $y=\min(l-k,(b-d)_+),H=k+a-e+(d-b)_+$.

\item[(ii)] If $a<e$ and $d-b\le e-a\le l-k$,
\begin{eqnarray*} 
R: z^m(f)\ot z^{m'}(d,c,b)\mapsto&&\\
&&\hspace{-3.8cm}z^{m'+H}(k-2a)\ot z^{m-H}(d+l-k-e+a-z,c,b+e-a-z),
\end{eqnarray*}
where $z=\min(b-d+e-a,l-k-e+a),H=k$.

\item[(iii)] Otherwise,
\begin{eqnarray*}
R: z^m(f)\ot z^{m'}(d,c,b)\mapsto&&\\
&&\hspace{-2.4cm}z^{m'+H}(k-2e+2w)\ot z^{m-H}(d+l-k-w,c,b+w),
\end{eqnarray*}
where $w=\min(l-k,(2e-2a-d+b)_+),H=k+\max(e-a-l+k,d-b-e+a)$.
\end{itemize}
\end{itemize}
\end{proposition}

Here and in what follows we use the notation
\[
(x)_{\pm} = \max(0,\pm x).
\]

\begin{lemma}[\cite{HKOT1,HKOT2}] \label{lem:deg0}
Under the isomorphism $B_l \ot B_1 \simeq B_1 \ot B_l$ we have
\begin{align*}
(l-j,j,0) \ot \hako{1} &\simeq 
\begin{cases}
\hako{2}\ot(l-j+1,j-1,0) & 0 < j \le l,\\
\hako{1}\ot(l,0,0) & j=0,
\end{cases}\\
(l-j,j,0) \ot \hako{2} &\simeq 
\begin{cases}
\hako{1}\ot(l-j-1,j+1,0) & 0 \le j < l,\\
\hako{2}\ot(0,l,0) & j=l.
\end{cases}
\end{align*}
\end{lemma}

This is valid for all Types I, II and III,
and it is essentially the same with the $A^{(1)}_n$ case.

We also give all the images of the isomorphism $B_\natural\ot B_1\longrightarrow
B_1\ot B_\natural$ that will be needed in the next subsection.
Let the symbol

\begin{center}
\vertex({\alpha\atop\beta},\gamma,\gamma',{\alpha'\atop\beta'})
\end{center}
signify ${\alpha\choose\beta}\ot\hako{\gamma}\mapsto\hakoc{\gamma'}\ot
{\alpha'\choose\beta'}$ under the isomorphism. If $i\neq j(i,j=1,2,\cd,n)$,
we easily obtain
\begin{center}
\vertex({i\atop\ol{j}},i,i,{i\atop\ol{j}})\quad
\vertex({i\atop\ol{j}},\ol{j},\ol{j},{i\atop\ol{j}})
\end{center}
{}from the weight reason. We list up the other cases below.
\bigskip

Common for Type I and II:

\vertex({1\atop2},3,1,{2\atop3})
\vertex({1\atop3},2,3,{1\atop2})
\vertex({2\atop3},1,2,{1\atop3})
\vertex({2\atop\ol{2}},1,2,{1\atop\ol{2}})
\vertex({2\atop\ol{1}},1,2,{2\atop\ol{2}})

\vertex({3\atop\ol{1}},1,3,{2\atop\ol{2}})
\vertex({3\atop\ol{1}},2,3,{2\atop\ol{1}})
\vertex({\ol{2}\atop\ol{1}},1,\ol{2},{2\atop\ol{2}})
\vertex({\ol{2}\atop\ol{1}},2,\ol{2},{2\atop\ol{1}})
\vertex({\ol{2}\atop\ol{1}},3,\ol{2},{3\atop\ol{1}})

For Type I only:

\vertex({1\atop2},\ol{2},1,\phi)
\vertex(\phi,3,1,{3\atop\ol{1}})
\vertex(\phi,\ol{2},1,{\ol{2}\atop\ol{1}})

For Type II only:

\vertex({1\atop2},\ol{2},1,{2\atop\ol{2}})
\vertex({1\atop2},\ol{1},1,{2\atop\ol{1}})
\vertex({1\atop\ol{2}},2,\ol{2},{1\atop2})
\vertex({1\atop\ol{2}},3,\ol{2},{1\atop3})

\vertex({1\atop\ol{2}},\ol{1},1,{\ol{2}\atop\ol{1}})
\vertex({2\atop\ol{2}},\ol{2},\ol{1},{1\atop\ol{2}})
\vertex({2\atop\ol{1}},3,\ol{1},{2\atop3})
\vertex({2\atop\ol{1}},\ol{2},\ol{1},{2\atop\ol{2}})

\subsection{Type I and II}
Let $p\in B_1^{\ot N}$ be a state. For an element $v\in B_\natural$, define
a map $:B_1^{\ot N}\longrightarrow B_1^{\ot N}$ by 
\begin{equation}\label{eq:natural}
\begin{array}{ccc}
B_\natural\ot B_1^{\ot N}&\stackrel{\sim}{\longrightarrow}&
B_1^{\ot N}\ot B_\natural\\
v\ot p&\mapsto&T_\natural[v](p)\ot v',
\end{array}
\end{equation}
where $v'$ is uniquely determined {}from $v$ and $p$ although it 
is not necessary for defining $T_\natural[v]$.

For Type I case, we have 
\begin{lemma} \label{lem:T-naturalI}
Consider the one-soliton state identified with $z^m(d,c,b)$ as an element of
$\Aff(B_k)$. Then we have
\[
T_\natural\left[{1\choose2}\right](z^m(d,c,b))=\left\{
\begin{array}{ll}
z^{m-2}(d+1,c,b-1)&\quad\mbox{if }b>0,\\
z^{m-1}(d+1,c-1,0)&\quad\mbox{if }b=0,c>0,\\
z^m(d,0,0)&\quad\mbox{if }b=c=0.
\end{array}\right.
\]
\end{lemma}

For Type II case, we have 

\begin{lemma} \label{lem:T-naturalII}
Consider the one-soliton state identified with $z^m(d,c,b)$ as an element of
$\Aff(B_k)$. Define $a$ by $k=2a+b+c+d$. Then we have

If $a=b=c=0$, then
\[
T_\natural\left[{1\choose2}\right](z^m(d,c,b))=z^m(k,0,0),
\]

else
\[
T_\natural\left[{1\choose2}\right](z^m(d,c,b))=\left\{
\begin{array}{ll}
z^{m-1}(d-1,c,b-1)&\quad\mbox{if }bd>0,\\
z^{m-1}(d+(b-1)_-,c-1,(b-1)_+)&\quad\mbox{if }bd=0,c>0,\\
z^{m-1}(d+(b-2)_-,0,(b-2)_+)&\quad\mbox{if }bd=c=0.
\end{array}\right.
\]
\[
T_\natural\left[{1\choose\ol{2}}\right](z^m(d,c,b))=\left\{
\begin{array}{ll}
z^{m-1}(d+1,c,b+1)&\quad\mbox{if }a>0,\\
z^m(d-1,c,b+1)&\quad\mbox{if }a=0,d>0,\\
z^m(0,c-1,b+1)&\quad\mbox{if }a=d=0,c>0,\\
z^m(0,0,b)&\quad\mbox{if }a=d=c=0.
\end{array}\right.
\]
\end{lemma}

Instead of proving these lemmas, we give examples below.

\begin{example}\label{ex:natural}
\begin{itemize}
\item[(i)] Type I ($k=4$)


\setlength{\unitlength}{1mm}
\begin{picture}(120,30)(0,0)
\multiput(8,15)(12,0){9}{\line(1,0){8}}
\multiput(12,9)(12,0){9}{\line(0,1){12}}

\put(5,14){${1\atop2}$}
\put(17,14){${1\atop2}$}
\put(29,14){${1\atop2}$}
\put(41,14){$\phi$}
\put(53,14){${\ol{2}\atop\ol{1}}$}
\put(65,14){${3\atop\ol{1}}$}
\put(77,14){${2\atop\ol{1}}$}
\put(89,14){${2\atop\ol{2}}$}
\put(101,14){${1\atop\ol{2}}$}
\put(113,14){${1\atop\ol{2}}$}

\put(11.4,22.3){$\scriptstyle 1$}
\put(23.4,22.3){$\scriptstyle 1$}
\put(35.4,22.3){$\scriptstyle \ol{2}$}
\put(47.4,22.3){$\scriptstyle \ol{2}$}
\put(59.4,22.3){$\scriptstyle 3$}
\put(71.4,22.3){$\scriptstyle 2$}
\put(83.4,22.3){$\scriptstyle 1$}
\put(95.4,22.3){$\scriptstyle 1$}
\put(107.4,22.3){$\scriptstyle 1$}

\put(11.4,5.9){$\scriptstyle 1$}
\put(23.4,5.9){$\scriptstyle 1$}
\put(35.4,5.9){$\scriptstyle 1$}
\put(47.4,5.9){$\scriptstyle 1$}
\put(59.4,5.9){$\scriptstyle \ol{2}$}
\put(71.4,5.9){$\scriptstyle 3$}
\put(83.4,5.9){$\scriptstyle 2$}
\put(95.4,5.9){$\scriptstyle 2$}
\put(107.4,5.9){$\scriptstyle 1$}
\end{picture}

Thus we have 
\[
T_\natural\left[{1\choose2}\right](z^m(1,1,2))=z^{m-2}(2,1,1).
\]

\item[(ii)] Type II ($k=6,a=1$)


\setlength{\unitlength}{1mm}
\begin{picture}(120,30)(0,0)
\multiput(8,15)(12,0){9}{\line(1,0){8}}
\multiput(12,9)(12,0){9}{\line(0,1){12}}

\put(5,14){${1\atop2}$}
\put(17,14){${1\atop2}$}
\put(29,14){${1\atop2}$}
\put(41,14){${2\atop\ol{1}}$}
\put(53,14){${2\atop\ol{2}}$}
\put(65,14){${1\atop\ol{2}}$}
\put(77,14){${1\atop\ol{2}}$}
\put(89,14){${1\atop{3}}$}
\put(101,14){${1\atop{3}}$}
\put(113,14){${1\atop{3}}$}

\put(11.4,22.3){$\scriptstyle 1$}
\put(23.4,22.3){$\scriptstyle 1$}
\put(35.4,22.3){$\scriptstyle \ol{1}$}
\put(47.4,22.3){$\scriptstyle \ol{2}$}
\put(59.4,22.3){$\scriptstyle \ol{2}$}
\put(71.4,22.3){$\scriptstyle \ol{2}$}
\put(83.4,22.3){$\scriptstyle 3$}
\put(95.4,22.3){$\scriptstyle 1$}
\put(107.4,22.3){$\scriptstyle 1$}

\put(11.4,5.9){$\scriptstyle 1$}
\put(23.4,5.9){$\scriptstyle 1$}
\put(35.4,5.9){$\scriptstyle 1$}
\put(47.4,5.9){$\scriptstyle \ol{1}$}
\put(59.4,5.9){$\scriptstyle \ol{1}$}
\put(71.4,5.9){$\scriptstyle \ol{2}$}
\put(83.4,5.9){$\scriptstyle \ol{2}$}
\put(95.4,5.9){$\scriptstyle 1$}
\put(107.4,5.9){$\scriptstyle 1$}
\end{picture}

Thus we have 
\[
T_\natural\left[{1\choose2}\right](z^m(0,1,3))=z^{m-1}(0,0,2).
\]

\item[(iii)] Type II ($k=6,a=2$)


\setlength{\unitlength}{1mm}
\begin{picture}(120,30)(0,0)
\multiput(8,15)(12,0){8}{\line(1,0){8}}
\multiput(12,9)(12,0){8}{\line(0,1){12}}

\put(5,14){${1\atop2}$}
\put(17,14){${1\atop2}$}
\put(29,14){${2\atop\ol{1}}$}
\put(41,14){${2\atop\ol{1}}$}
\put(53,14){${2\atop\ol{1}}$}
\put(65,14){${2\atop\ol{1}}$}
\put(77,14){${2\atop\ol{2}}$}
\put(89,14){${1\atop\ol{2}}$}
\put(101,14){${1\atop\ol{2}}$}

\put(11.4,22.3){$\scriptstyle 1$}
\put(23.4,22.3){$\scriptstyle \ol{1}$}
\put(35.4,22.3){$\scriptstyle \ol{1}$}
\put(47.4,22.3){$\scriptstyle 2$}
\put(59.4,22.3){$\scriptstyle 2$}
\put(71.4,22.3){$\scriptstyle 1$}
\put(83.4,22.3){$\scriptstyle 1$}
\put(95.4,22.3){$\scriptstyle 1$}

\put(11.4,5.9){$\scriptstyle 1$}
\put(23.4,5.9){$\scriptstyle 1$}
\put(35.4,5.9){$\scriptstyle \ol{1}$}
\put(47.4,5.9){$\scriptstyle 2$}
\put(59.4,5.9){$\scriptstyle 2$}
\put(71.4,5.9){$\scriptstyle 2$}
\put(83.4,5.9){$\scriptstyle 2$}
\put(95.4,5.9){$\scriptstyle 1$}

\end{picture}

Thus we have 
\[
T_\natural\left[{1\choose2}\right](z^m(2,0,0))=z^{m-1}(4,0,0).
\]
\end{itemize}
\end{example}

In what follows, we understand 
$T_\natural=T_\natural\left[{1\choose2}\right]$.
Let 
\begin{align*}
p_1 &= .....\imath_{k_1}(b_1)......................................\\
p_2 &= .....\imath_{k_1}(b_1).............\imath_{k_2}(b_2)...............
\end{align*}
be the one and two-soliton states which are identified with 
$z^{m_1}b_1 \in \Aff(B_{k_1})$ and 
$z^{m_1}b_1 \ot z^{m_2}b_2 \in \Aff(B_{k_1}) \ot \Aff(B_{k_2})$,
respectively.
Suppose that ${1 \choose 2} \ot p_1 \simeq T_\natural(p_1) \ot v'$
as in (\ref{eq:natural}).
Then $T_\natural(p_1)$ is another one-soliton state 
identified with some element $T_\natural(z^{m_1}b_1) \in \Aff(B_{k_1})$, 
which is explicitly given in {\sc Lemmas} \ref{lem:T-naturalI} and \ref{lem:T-naturalII} 
for $b_1$ of the form (\ref{eq:notation}).
Similarly, 
$T_\natural(p_2)$ is a two-soliton state 
identified with 
$T_\natural(z^{m_1}b_1)\ot T_\natural[v'](z^{m_2}b_2) 
\in \Aff(B_{k_1}) \ot  \Aff(B_{k_2})$.
When $b_1 \ot b_2$ is a highest weight element,
we will only  encounter  $v'= {1 \choose 2}$ or ${1 \choose \ol{2}}$.
In fact one always has $v' = {1 \choose 2}$ for Type I (and $A^{(1)}_n$).
On the other hand for Type II,
$v' = {1 \choose \ol{2}}$ happens if and only if 
$b_1 = (f)$ with $f \le k_1-2$, and 
$v' = {1 \choose 2}$ otherwise.
See {\sc Proposition} \ref{prop:hweII} (1) for the form of the
highest weight elements and 
{\sc Example} \ref{ex:natural} (iii) for a situation 
$v' = {1 \choose \ol{2}}$.
Thus the results on $T_\natural[{1\choose 2}]$ and 
$T_\natural[{1\choose \ol{2}}]$ 
in {\sc Lemmas} \ref{lem:T-naturalI} and \ref{lem:T-naturalII}
are enough to calculate the action of $T_\natural[{1 \choose 2}]$
on any highest two-soliton state $z^{m_1}b_1 \ot z^{m_2}b_2$.
Later in the proof of  {\sc Lemma} \ref{lem:Comm-R-T-natural},
we will extensively use Lemmas \ref{lem:T-naturalI} and \ref{lem:T-naturalII}
either  for $(k_1,k_2) = (l,k)$ or $(k,l)$.

Let us prepare a few more facts.
\begin{proposition}\label{pr:Tcommute}
\[
T_\natural T_r = T_r T_\natural\quad \text{for any } \, r \ge 1.
\]
\end{proposition}
\Proof
By the weight reason we have
${1 \choose 2 } \ot u_r \simeq u_r \ot {1 \choose 2 }$ under the 
isomorphism $B_{\natural} \ot B_r \simeq B_r \ot B_\natural$.
\qed
\begin{lemma} \label{lem:Comm-R-T-natural}
Assume $l>k$.
Identify a two-soliton state with an element of $\Aff(B_l)\ot\Aff(B_k)$.
Then on the set of highest weight elements of $\Aff(B_l)\ot\Aff(B_k)$,
the map $T_\natural$ commutes with the combinatorial $R$.
\end{lemma}

\Proof 
Consider Type I case first. 
We apply $T_\natural$ on both sides of (\ref{eq:R-hwe}). In view of 
{\sc Lemma} \ref{lem:T-naturalI}, we divide the check into 3 cases.

\begin{itemize}
\item[(i)] $b>0$:
\begin{eqnarray*}
LHS&=&z^m(l)\ot z^{m'-2}(d+1,c,b-1),\\
RHS&=&z^{m'+H}(k)\ot z^{m-H-2}(l-b-c+1,c,b-1).
\end{eqnarray*}

\item[(ii)] $b=0,c>0$:
\begin{eqnarray*}
LHS&=&z^m(l)\ot z^{m'-1}(d+1,c-1,0),\\
RHS&=&z^{m'+H}(k)\ot z^{m-H-1}(l-c+1,c-1,0).
\end{eqnarray*}

\item[(iii)] $b=c=0$: trivial.

\end{itemize}
By {\sc Proposition} \ref{prop:hweI} (2), we see the LHS is sent to the RHS by 
the combinatorial $R$.

The proof of Type II case gets more involved.
Set $p=z^m(f)\ot z^{m'}(d,c,b)$ and 
$p'=R(p)=z^{m'+H}(\varphi)\ot z^{m-H}(\delta,\gamma,\beta)$. Recall integers
$a,e$ are defined by $l=2e+f,k=2a+b+c+d$. First consider the case when
$\min(a,e)>0$. Then we have

\begin{equation} \label{eq:A}
T_\natural(p)=z^{m-1}(f+2)\ot z^{m'-1}(d+1,c,b+1).
\end{equation}
On the other hand, in all cases (i)(ii)(iii) in 
{\sc Proposition} \ref{prop:hweII} (2) one finds $\varphi<k,
\beta+\gamma+\delta<l$. (Notice that $w<e-a$ in (iii).) Thus we have

\begin{equation} \label{eq:B}
T_\natural(p')=z^{m'+H-1}(\varphi+2)\ot z^{m-H-1}(\delta+1,\gamma,\beta+1).
\end{equation}
One can check (\ref{eq:A}) is sent to (\ref{eq:B}) by the combinatorial $R$
using {\sc Proposition} \ref{prop:hweII} (2) again. (Notice that 
$T_\natural(p)$ belongs to the same case as $p$.)

Now assume $\min(a,e)=0$. We divide the cases according to those in 
{\sc Proposition} \ref{prop:hweII} (2). In all cases we list $T_\natural(p)$
and $T_\natural(p')$ in the manner $T_\natural(p)\mapsto T_\natural(p')$.
Checking $T_\natural(p)$ is sent to $T_\natural(p')$ by $R$ is left to the
reader.

\begin{itemize}
\item[(i)] $a\ge e$;
We have $e=0$, and $p'=z^{m'+H}(k)\ot z^{m-H}(\delta,c,\beta)$ with
$\delta=d+l-k-y,\beta=b-y,y=\min(l-k,(b-d)_+),H=k+a+(d-b)_+$.

\begin{itemize}
\item[$\bullet$] $bd>0$;
The condition implies $\beta\delta>0$.
\[
z^m(l)\ot z^{m'-1}(d-1,c,b-1)\mapsto 
z^{m'+H}(k)\ot z^{m-H-1}(\delta-1,c,\beta-1).
\]

\item[$\bullet$] $b=0$;
We have $y=0,\delta=d+l-k,\beta=0$.

If $c>0$,
\[
z^m(l)\ot z^{m'-1}(d+1,c-1,0)\mapsto 
z^{m'+H}(k)\ot z^{m-H-1}(\delta+1,c-1,0).
\]

If $c=0, d \le k-2$,
\[
z^m(l)\ot z^{m'-1}(d+2,0,0)\mapsto 
z^{m'+H}(k)\ot z^{m-H-1}(\delta+2,0,0).
\]

If $c=0, d = k$,
\[
z^m(l)\ot z^{m'}(k) \mapsto
z^{m'+H}(k)\ot z^{m-H}(l).
\]

\item[$\bullet$] $b>0,d=0,l-k<b$;
We have $y=l-k,\delta=0,\beta=b-l+k(\ge1)$. Note $b\ge2$.

If $c>0$,
\[
z^m(l)\ot z^{m'-1}(0,c-1,b-1)\mapsto 
z^{m'+H}(k)\ot z^{m-H-1}(0,c-1,\beta-1).
\]

If $c=0,\beta\ge2$,
\[
z^m(l)\ot z^{m'-1}(0,0,b-2)\mapsto 
z^{m'+H}(k)\ot z^{m-H-1}(0,0,\beta-2).
\]

If $c=0,\beta=1$,
\[
z^m(l)\ot z^{m'-1}(0,0,b-2)\mapsto 
z^{m'+H}(k)\ot z^{m-H-1}(1,0,0).
\]

\item[$\bullet$] $b>0,d=0,l-k\ge b$;
We have $y=b,\delta=l-k-b,\beta=0$.

If $c>0$,
\[
z^m(l)\ot z^{m'-1}(0,c-1,b-1)\mapsto 
z^{m'+H}(k)\ot z^{m-H-1}(\delta+1,c-1,0).
\]

If $c=0,b\ge2$,
\[
z^m(l)\ot z^{m'-1}(0,0,b-2)\mapsto 
z^{m'+H}(k)\ot z^{m-H-1}(\delta+2,0,0).
\]

If $c=0,b=1$,
\[
z^m(l)\ot z^{m'-1}(1,0,0)\mapsto 
z^{m'+H}(k)\ot z^{m-H-1}(\delta+2,0,0).
\]
\end{itemize}

\item[(ii)] $a<e,l-k\ge e-a\ge d-b$; 
We have $a=0$, and $p'=z^{m'+H}(k)\ot z^{m-H}(\delta,c,\beta)$ with
$\delta=d+l-k-e-z,\beta=b+e-z,z=\min(b-d+e,l-k-e),H=k$. 
If $b-d+e\le l-k-e$, then $z=b-d+e,\delta\ge\beta=d$, else, 
$z=l-k-e,\beta>\delta=d$.

\begin{itemize}
\item[$\bullet$] $d>0$;
\[
z^{m-1}(f+2)\ot z^{m'}(d-1,c,b+1)\mapsto 
z^{m'+H}(k)\ot z^{m-H-1}(\delta-1,c,\beta-1).
\]

\item[$\bullet$] $d=0,c>0$;

If $b-d+e\le l-k-e$,
\[
z^{m-1}(f+2)\ot z^{m'}(0,c-1,b+1)\mapsto 
z^{m'+H}(k)\ot z^{m-H-1}(\delta+1,c-1,0).
\]

If $b-d+e>l-k-e$,
\[
z^{m-1}(f+2)\ot z^{m'}(0,c-1,b+1)\mapsto 
z^{m'+H}(k)\ot z^{m-H-1}(0,c-1,\beta-1).
\]

\item[$\bullet$] $d=c=0$;

If $b-d+e\le l-k-e$,
\[
z^{m-1}(f+2)\ot z^{m'}(0,0,b)\mapsto 
z^{m'+H}(k)\ot z^{m-H-1}(\delta+2,0,0).
\]

If $b-d+e>l-k-e,\beta\ge2$,
\[
z^{m-1}(f+2)\ot z^{m'}(0,0,b)\mapsto 
z^{m'+H}(k)\ot z^{m-H-1}(0,0,\beta-2).
\]

If $b-d+e>l-k-e,\beta=1$,
\[
z^{m-1}(f+2)\ot z^{m'}(0,0,b)\mapsto 
z^{m'+H}(k)\ot z^{m-H-1}(1,0,0).
\]
\end{itemize}

\item[(iii)] $a<e,l-k<e-a\mbox{ or }e-a<d-b$;
We have $a=0$, and $p'=z^{m'+H}(\varphi)\ot z^{m-H}(\delta,c,\beta)$ with
$\varphi=k-2e+2w,\delta=d+l-k-w,\beta=b+w,w=\min(l-k,(2e-d+b)_+),
H=k+\max(e-l+k,d-b-e)$. Note that $w<e$, and if $d=0$, $w=l-k$.

\begin{itemize}
\item[$\bullet$] $d>0$;
\[
z^{m-1}(f+2)\ot z^{m'}(d-1,c,b+1)\mapsto 
z^{m'+H-1}(\varphi+2)\ot z^{m-H}(\delta-1,c,\beta+1).
\]

\item[$\bullet$] $d=0,c>0$;
\[
z^{m-1}(f+2)\ot z^{m'}(0,c-1,b+1)\mapsto 
z^{m'+H-1}(\varphi+2)\ot z^{m-H}(0,c-1,\beta+1).
\]

\item[$\bullet$] $d=c=0$;
\[
z^{m-1}(f+2)\ot z^{m'}(0,0,b)\mapsto 
z^{m'+H-1}(\varphi+2)\ot z^{m-H}(0,0,\beta).
\]
\end{itemize}
\end{itemize}
\qed

\begin{lemma} \label{lem:wt-difference}
For arbitrary  $l,k \in \Z_{\ge 1}$  
(not necessarily $l > k$ as assumed previously),
let $z^{m_i}b_i\ot z^{m'_i}b'_i$ ($i=1,2$) be highest weight elements in
$\Aff(B_l)\ot\Aff(B_k)$. If $T_\natural(z^{m_1}b_1\ot z^{m'_1}b'_1)=
T_\natural(z^{m_2}b_2\ot z^{m'_2}b'_2)$ and 
$\wt(b_1\ot b'_1)=\wt(b_2\ot b'_2)$, then $z^{m_1}b_1\ot z^{m'_1}b'_1=
z^{m_2}b_2\ot z^{m'_2}b'_2$.
\end{lemma}

\Proof 
We first consider Type I case. 
{}From {\sc Lemma} \ref{lem:T-naturalI} we see the inverse image of an element
of the form $z^m(l)\ot z^{m'}(d,c,b)$ has more than one element only
when $d>0,b=0$. In this case $T_\natural^{-1}(z^m(l)\ot z^{m'}(d,c,b))$ 
consists of the following elements:
\begin{eqnarray*}
&&z^m(l)\ot z^{m'+2}(d-1,c,1),\\
&&z^m(l)\ot z^{m'+1}(d-1,c+1,0),\\
&&z^m(l)\ot z^{m'}(d,0,0)\quad\mbox{if }c=0.
\end{eqnarray*}
It is easy to see the weights are distinct.

For Type II the cases when $T_\natural^{-1}(z^m(f)\ot z^{m'}(d,c,b))$ has
more than one element are the following. Using {\sc Lemma} \ref{lem:T-naturalII}
we divide into 6 cases.
\begin{itemize}
\item[(i)] $d>0,b>0,a>0,e=0$:
\begin{eqnarray*}
&&z^m(l)\ot z^{m'+1}(d+1,c,b+1),\\
&&z^{m+1}(l-2)\ot z^{m'+1}(d-1,c,b-1).
\end{eqnarray*}

\item[(ii)] $d>0,b>0,a=0$:
\begin{eqnarray*}
&&z^{m+1}(f-2)\ot z^{m'+1}(d-1,c,b-1),\\
&&z^{m+1}(f-2)\ot z^{m'}(d+1,c,b-1).
\end{eqnarray*}

\item[(iii)] $d>0,b=0,c>0,a>0,e=0$:
\begin{eqnarray*}
&&z^m(l)\ot z^{m'+1}(d+1,c,1),\\
&&z^m(l)\ot z^{m'+1}(d-1,c+1,0).
\end{eqnarray*}

\item[(iv)] $d>0,b=0,c=0,e=0$:
\begin{eqnarray*}
&&z^m(l)\ot z^{m'+1}(d-1,1,0),\\
&&z^m(l)\ot z^{m'+1}(d-2,0,0)\quad\mbox{if }d\ge2,\\
&&z^m(l)\ot z^{m'+1}(0,0,1)\quad\mbox{if }d=1,\\
&&z^m(l)\ot z^{m'+1}(d+1,0,1)\quad\mbox{if }a>0,\\
&&z^m(l)\ot z^{m'}(d,0,0)\quad\mbox{if }a=0.
\end{eqnarray*}

\item[(v)] $d=0,b>0,a=0$:
\begin{eqnarray*}
&&z^{m+1}(f-2)\ot z^{m'}(0,c+1,b-1),\\
&&z^{m+1}(f-2)\ot z^{m'}(1,c,b-1),\\
&&z^{m+1}(f-2)\ot z^{m'}(0,0,b)\quad\mbox{if }c=0.
\end{eqnarray*}

\item[(vi)] $d=0,a>0,e=0$:
\begin{eqnarray*}
&&z^m(l)\ot z^{m'+1}(1,c,b+1),\\
&&z^m(l)\ot z^{m'+1}(0,c+1,b+1),\\
&&z^m(l)\ot z^{m'+1}(0,0,b+2)\quad\mbox{if }c=0.
\end{eqnarray*}

\end{itemize}
In all cases the weights are distinct.
\qed
\vskip0.3cm
 
We introduce a degree on the set of highest weight elements $\Aff(B_l)\ot\Aff(B_k)$. 
First, consider an element $p$ such that $T_\natural(p)=p$. 
We define the degree of such $p$,
denoted by $\deg p$, to be 0. {}From {\sc Lemmas} \ref{lem:T-naturalI} 
and  \ref{lem:T-naturalII} we easily see that 
two-soliton states of degree 0 must be of the form $z^m(l)\ot z^{m'}(k)$
for some $m$ and $m'$.
{}From {\sc Propositions} \ref{prop:hweI} and \ref{prop:hweII} we also see that 
any other highest weight elements can be 
transformed to a degree 0 one by applying $T_\natural$ for 
appropriate  times.
Thus by setting $\deg T_\natural(p)=\deg p-1$, 
we can assign a non-negative degree uniquely to all the highest weight 
elements.

\vskip0.2cm
{\it Proof of {\sc Theorem} \ref{th:main} for Type I and II.}
Identify a two-soliton state $p$ with an element of 
$\Aff(B_l)\ot\Aff(B_k)$.
For $r(>k)$ and sufficiently large $t$  we are to show 
$T_r^t(p)=R(p)$.
Due to the $U_q(\gehb_{n-1})$ invariance mentioned in 
Section \ref{subsection:pre},
it suffices to check the assertion when $p$ is a highest weight element.
We do this by induction on the degree. 
If $\deg p=0$, the statement can be directly checked by means of 
{\sc Lemma} \ref{lem:deg0}. 
In fact the calculation in this case is the same as 
$A_n^{(1)}$ considered in Section \ref{subsection:Atype}. 
Now assume $\deg p>0$. We have 
\[
(T_\natural T_r^t)(p)\stackrel{\rm {\sc Prop}\;\ref{pr:Tcommute}}{=}
(T_r^t T_\natural)(p)=
(RT_\natural)(p)\stackrel{\rm {\sc Lemma}\;\ref{lem:Comm-R-T-natural}}{=}
(T_\natural R)(p).
\]
Here the second equality is assured by the induction hypothesis.
Since both $T_r^t(p)$ and $R(p)$  have the same weight as $p$, 
the proof is finished by {\sc Lemma} \ref{lem:wt-difference}.
\qed

\subsection{\mathversion{bold}Type III : $A_{2n}^{(2)}$ and $D^{(2)}_{n+1}$ cases}
\label{subsection:typeIII}

Let us proceed to the case 
$\geh = A^{(2)}_{2n},D^{(2)}_{n+1}$.
In this subsection $B_l$ denotes the $U'_q(\geh)$-crystal 
and  ${\tilde B}_l$ does the $U'_q(C^{(1)}_n)$-crystal.

First we quote some results {}from \cite{HKOT2}.
For $\geh = A^{(2)}_{2n},D^{(2)}_{n+1}$ and for 
any $l \in {\mathbb Z}_{\ge 1}$, let  $\omega$ be the map
defined by
\begin{align*}
\omega: \; B_l \; &\longrightarrow \;
{\tilde B}_{2l}\\
(x_1,\ldots,x_n,&\ol{x}_n,\ldots,\ol{x}_1)\\
 &\mapsto (2x_1,\ldots,2x_n,2\ol{x}_n,\ldots,2\ol{x}_1)
\;\; \text{for } \geh = A^{(2)}_{2n},\\
(x_1,\ldots,x_n,&x_0,\ol{x}_n,\ldots,\ol{x}_1)\\
 &\mapsto (2x_1,\ldots,2x_{n-1},2x_n+x_0,
 2\ol{x}_n+x_0,2\ol{x}_{n-1},\ldots,2\ol{x}_1)
\; \; \text{for } \geh = D^{(2)}_{n+1}.
\end{align*}
We do not exhibit the $l$-dependence in $\omega$.
The map $\omega$ for $\geh_{n-1}$ will be denoted as $\omega_{n-1}$.

\begin{theorem}[\cite{HKOT2}]\label{th:RforAD}
Let $\geh = A^{(2)}_{2n}, D^{(2)}_{n+1}$.
Suppose $b \ot c \stackrel{\sim}{\mapsto} {\tilde c} \ot {\tilde b}$
under the isomorphism $B_l \ot B_k 
\stackrel{\sim}{\rightarrow} B_k \ot B_l$
of $U'_q(\geh)$-crystals.
Then one has 
\begin{itemize}
\item[(1)]
$\omega(b) \ot \omega(c)\stackrel{\sim}{\mapsto} 
\omega({\tilde c}) \ot \omega({\tilde b})$
under the isomorphism ${\tilde B}_{2l} \ot {\tilde B}_{2k} 
\stackrel{\sim}{\rightarrow} {\tilde B}_{2k} \ot {\tilde B}_{2l}$
of $U'_q(C^{(1)}_n)$-crystals.
\item[(2)]
$H^{\geh}(b\ot c) = H^{C^{(1)}_n}(\omega(b) \ot \omega(c))$.
\end{itemize}
\end{theorem}

For $\geh = A^{(2)}_{2n}$ and $D^{(2)}_{n+1}$ we introduce a map $\eta$ as

\begin{align*}
\eta: &\quad B_1 \quad \;\; \longrightarrow \quad
{\tilde B}_1 \ot {\tilde B}_1\\
&\quad \hako{a} \quad\quad\; \mapsto \qquad \hako{a} \ot \hako{a}\quad 
(a \in \{1,\ldots,n,\ol{n},\ldots, \ol{1}\})\\
&\quad\; \phi \quad\quad\;\; \mapsto \qquad \hako{1} \ot \hako{\ol{1}}\\
&\quad \hako{0} \quad\quad\; \mapsto \qquad \hako{\ol{n}} \ot \hako{n} 
\quad (\mbox{ only for } \geh = D^{(2)}_{n+1}).
\end{align*}

\begin{lemma}\label{lem:fact2}
For $\geh = A^{(2)}_{2n}$ and $D^{(2)}_{n+1}$ the following 
diagram is commutative for any $L \in {\mathbb Z}_{\ge 1}$.
\begin{equation*}
\begin{CD}
B_L \ot B_1 @>{\sim}>> B_1 \ot B_L
\\
@V{\omega \ot \eta}VV  @VV{\eta \ot \omega}V \\
{\tilde B}_{2L} \ot {\tilde B}_{1}\ot {\tilde B}_{1}
 @>{\sim}>>
{\tilde B}_{1}\ot {\tilde B}_{1} \ot {\tilde B}_{2L}
\end{CD}
\end{equation*}
\end{lemma}
\Proof
For $\geh = A^{(2)}_{2n}$ (resp. $D^{(2)}_{n+1}$) 
the image $\eta(B_1) \subset {\tilde B}_{1}\ot {\tilde B}_{1}$ 
is stable under the actions of 
$\et{i}^2, \ft{i}^2 (i \neq 0), \et{0}$ and $\ft{0}$
(resp. $\et{i}^2, \ft{i}^2 (i \neq 0,n), \et{0}, \et{n}, \ft{0}$ and $\ft{n}$).
Moreover these operators generate the identical crystal subgraph within 
${\tilde B}_{1}\ot {\tilde B}_{1}$
as they do in $\omega(B_1) \subset {\tilde B}_2$.
The assertion follows {}from this fact and {\sc Theorem} \ref{th:RforAD}.
\qed

For $\geh = A^{(2)}_{2n}, D^{(2)}_{n+1}$ case, denote the map  
$\imath_l: B_l \rightarrow \left(B_1\right)^{\ot l}$
specified in Section \ref{subsection:solitons} by $\imath_l$, 
and for $C^{(1)}_n$ case by ${\tilde \imath}_l$.
Then the following can be checked by a direct calculation.
\begin{lemma}\label{lem:fact3}
Let $\geh = A^{(2)}_{2n}$ or $D^{(2)}_{n+1}$, and $\hako{1} \in {\tilde B}_1$.
For any $b \in B_l$ for $U'_q(\geh_{n-1})$ we have
\begin{equation*}
\eta^{\ot l}\left(\imath_l(b)\right) \ot \hako{1} =
\begin{cases}
{\tilde \imath}_{2l}(\omega_{n-1}(b)) \ot \hako{1} & \mbox{if $l-s(b)$ is even},\\
\hako{1}\ot {\tilde \imath}_{2l}(\omega_{n-1}(b))& \mbox{if $l-s(b)$ is odd},
\end{cases}
\end{equation*}
where $s(b)$ has been defined in Section \ref{subsection:solitons}.
\end{lemma}

{\it Proof of {\sc Theorem} \ref{th:main} for Type III}.
We attribute the proof for
$\geh = A^{(2)}_{2n}, D^{(2)}_{n+1}$ 
to the $C^{(1)}_n$ case treated in the 
previous subsection.
Suppose that
\begin{equation}\label{eq:b1b2}
z^{\xi_1}b_1 \ot z^{\xi_2}b_2 \stackrel{\sim}{\mapsto}
z^{\xi_2+h}\tilde{b}_2 \ot z^{\xi_1-h}\tilde{b}_1 
\end{equation}
under the isomorphism 
$\Aff(B_{l})\ot \Aff(B_{k})
\stackrel{\sim}{\rightarrow}
\Aff(B_{k})\ot \Aff(B_{l})$
of $U_q(\geh_{n-1})$-crystals.
Then {}from {\sc Theorem} \ref{th:RforAD} we have
\begin{equation}\label{eq:cb1b2}
z^{\xi_1}\omega_{n-1}(b_1) \ot z^{\xi_2}\omega_{n-1}(b_2) \stackrel{\sim}{\mapsto}
z^{\xi_2+h}\omega_{n-1}(\tilde{b}_2) \ot z^{\xi_1-h}\omega_{n-1}(\tilde{b}_1)
\end{equation}
under the isomorphism 
$\Aff({\tilde B}_{2l})\ot \Aff({\tilde B}_{2k})
\stackrel{\sim}{\rightarrow}
\Aff({\tilde B}_{2k})\ot \Aff({\tilde B}_{2l})$
of $U_q(C^{(1)}_{n-1})$-crystals.
Define $e,e',a,a'$ by
\begin{equation*}
e = l-s(b_1),\quad e' = l - s(\tilde{b}_1),\quad
a = k-s(b_2),\quad a' = k - s(\tilde{b}_2),
\end{equation*}
where $s$ for $\geh_{n-1}$ has been defined in Section \ref{subsection:solitons}.
Put
\begin{equation*}
p_1 =  \cd \ot 1 \ot \imath_l(b_1)_{m_1}\ot 1 \ot \cd \ot 1
\ot \imath_k(b_2)_{m_2}\ot 1 \ot \cd 
\in \left(B_1\right)^{\ot N},
\end{equation*}
where the notation $\imath_l(b_1)_{m_1}$ for example signifies that 
$\nu_1$ appearing in 
$\imath_l(b_1) = \nu_1 \ot \cd \ot \nu_l$ is located 
on the $m_1$-th component counting {}from the right end in $p_1$.
We assume that $N-m_1, m_1-m_2, m_2 \gg l>k$.
Let $t$ be sufficiently large and $L > k$.
Then {}from {\sc Lemmas} \ref{lem:1}, \ref{lem:2} and \ref{lem:fact2}
we have the commutative diagram:
\begin{equation*}
\begin{CD}
u^{\ot t}_L \ot p_1 @>{\sim}>> \tilde{p}_1 \ot u^{\ot t}_L
\\
@V{\omega^{\ot t} \ot \eta^{\ot N}}VV  
@VV{\eta^{\ot N} \ot \omega^{\ot t}}V \\
u^{\ot t}_{2L} \ot p_2 @>{\sim}>> \tilde{p}_2 \ot u^{\ot t}_{2L}
\end{CD}
\end{equation*}
Here the upper (lower) $\stackrel{\sim}{\rightarrow}$ denotes 
the isomorphism of $U'_q(\geh) (U'_q(C^{(1)}_n))$-crystals.
According to (\ref{eq:coef2-1})--(\ref{eq:coef2-2})
the incoming state $p_1$ corresponds to 
\begin{equation*}
z^{2m_1-\bar{e}}b_1 \ot z^{2m_2-\bar{a}}b_2
\in \Aff(B_{l})\ot \Aff(B_{k}),
\end{equation*}
where the RHS is a tensor product of $U_q(\geh_{n-1})$-crystal.
$\bar{a} = (1-(-1)^a)/2$ and similarly for $\bar{e}$.
Thus {}from (\ref{eq:b1b2}) we are to show 
that the outgoing state $\tilde{p}_1$ corresponds to 
\begin{equation}\label{eq:out}
z^{2m_2-\bar{a}+h}\tilde{b}_2 \ot z^{2m_1-\bar{e}-h}\tilde{b}_1
\in \Aff(B_{k})\ot \Aff(B_{l}),
\end{equation}
To do this we utilize the commutativity of the diagram
to explicitly determine $\tilde{p}_1$ as follows.
First we use {\sc Lemma} \ref{lem:fact3} to find 
\begin{align}
&p_2 = \cd \ot 1 \ot 
{\tilde \imath}_{2l}(\omega_{n-1}(b_1))_{\mu_1}\ot 1 \ot \cd \ot 1 
\ot {\tilde \imath}_{2k}(\omega_{n-1}(b_2))_{\mu_2}
\ot 1 \ot \cd 
\in \left({\tilde B}_1\right)^{\ot 2N},\nonumber \\
&\mu_1 = 2m_1 - \bar{e}, \qquad \mu_2 = 2m_2 - \bar{a}.
\label{eq:mu}
\end{align}
Next we invoke {\sc Theorem} \ref{th:main} for $C^{(1)}_n$ case to derive
\begin{align}
&\tilde{p}_2 = \cd \ot 1 \ot 
{\tilde \imath}_{2k}(\omega_{n-1}(\tilde{b}_2))_{\tilde{\mu}_2}
\ot 1 \ot \cd \ot 1
\ot {\tilde \imath}_{2l}(\omega_{n-1}(\tilde{b}_1))_{\tilde{\mu}_1}
\ot 1 \ot \cd 
\in \left({\tilde B}_1\right)^{\ot 2N},\nonumber \\
&\tilde{\mu}_1 = \mu_1 - 2\min(L,l)t-h, \qquad \tilde{\mu}_2 = \mu_2 - 2kt+h,
\label{eq:mutilde}
\end{align}
where the phase shift is $h$ 
in (\ref{eq:b1b2}) because of {\sc Theorem} \ref{th:RforAD} (2).
Finally by using {\sc Lemma} \ref{lem:fact3} again we obtain
\begin{align}
&\tilde{p}_1 =  \cd \ot 1 \ot 
\imath_l(\tilde{b}_2)_{\tilde{m}_2}\ot 1 \ot \cd \ot 1
\ot \imath_k(\tilde{b}_1)_{\tilde{m}_1}\ot 1 \ot \cd 
\in \left(B_1\right)^{\ot N},\nonumber \\
&\tilde{m}_1 = \frac{\tilde{\mu}_1+\bar{e'}}{2},\qquad 
\tilde{m}_2 = \frac{\tilde{\mu}_2+\bar{a'}}{2}.\label{eq:mtilde}
\end{align}
According to (\ref{eq:coef2-1}) and (\ref{eq:coef2-2})
this corresponds to 
\begin{equation*}
z^{2(\tilde{m}_2+kt)-\bar{a'}} \tilde{b}_2 \ot
z^{2(\tilde{m}_1+\min(L,l)t)-\bar{e'}} \tilde{b}_1.
\end{equation*}
Combining (\ref{eq:mu}),(\ref{eq:mutilde}) and (\ref{eq:mtilde}),
we find that this is indeed equal to (\ref{eq:out}).
\qed

\section{Examples}\label{section:examples}

Let us present examples of soliton scattering.
Given a state $p$ of an automaton,
to find the time evolution $T_r(p)$  involves a series of calculations of
$B_r \ot B_1 \simeq B_1 \ot B_r$.
For $r$ finite this is a tedious task in general
but can be done according to the insertion scheme in \cite{HKOT1,HKOT2}.
The presented examples are actually so generated by a computer
by setting $r=12$.
However if  $r$ is large enough so that $T = T_r$,
the result in \cite{HKT2} enables one 
to compute $T(p)$ easily in an analogous manner to the $A^{(1)}_n$ case
where the simple algorithm (box-ball rule) is known \cite{TS,T,TNS}.
In all the examples below $T = T_{12}$ is valid.

First we give examples of two-soliton scattering
described by
$\Aff(B_5) \otimes \Aff(B_3)
\xrightarrow{\sim} \Aff(B_3) \otimes \Aff(B_5)$.

\begin{example}\label{ex:A(1)}
$\mathfrak{g}=A^{(1)}_3.$\\
\noindent%
$0:\cdots 1 1 \kakomiu{4 4 3 2 2} 1 1 1 1 1 1
 \kakomid{4 3 3} 1 1 1 1 1 1 1 1 1 1 1 1 1 1 1
1 1 1 1 1 1 1 1 1 1 1 1 1 1 1 1 \cdots$\\[\sukima]%
$1:\cdots 1 1 1 1 1 1 1 \kakomiu{4 4 3 2 2} 1 1 1 1
 \kakomid{4 3 3} 1 1 1 1 1 1 1 1 1 1 1 1
1 1 1 1 1 1 1 1 1 1 1 1 1 1 1 1 \cdots$\\[\sukima]%
$2:\cdots 1 1 1 1 1 1 1 1 1 1 1 1 \kakomiu{4 4 3 2 2}
 1 1 \kakomid{4 3 3} 1 1 1 1 1 1 1 1 1
1 1 1 1 1 1 1 1 1 1 1 1 1 1 1 1 \cdots$\\[\sukima]%
$3:\cdots 1 1 1 1 1 1 1 1 1 1 1 1 1 1 1 1 1
 \kakomiu{4 4 3 2 2} \kakomid{4 3 3} 1 1 1 1 1 1
 1 1 1 1 1 1 1 1 1 1 1 1 1 1 1 1 \cdots$\\[\sukima]%
$4:\cdots 1 1 1 1 1 1 1 1 1 1 1 1 1 1 1 1 1 1 1 1 1 1
 \kakomiu{3 2 2} \kakomiud{4 4} \kakomid{4} 3 3 1
 1 1 1 1 1 1 1 1 1 1 1 1 1 1 1 1 \cdots$\\[\sukima]%
$5:\cdots 1 1 1 1 1 1 1 1 1 1 1 1 1 1 1 1 1 1 1 1 1 1 1 1 1 3 2
 \kakomiu{2} \kakomiud{1 1 4} \kakomiu{4}
 4 3 3 1 1 1 1 1 1 1 1 1 1 1 1 \cdots$\\[\sukima]%
$6:\cdots 1 1 1 1 1 1 1 1 1 1 1 1 1 1 1 1 1 1 1 1 1 1 1 1 1 1 1 1 3 2 2
\kakomid{1} \kakomiud{1 1} \kakomiu{1 4 4} 4 3 3 1 1 1 1 1 1 1 \cdots$\\[\sukima]%

Here the markers specify the positions of initial solitons
under the interaction-free propagation.
One can read off the phase shift {}from the deviation of the
outgoing solitons {}from the corresponding markers.
In terms of the soliton labels, the above scattering  is described by
the combinatorial $R$  of $U_q(A^{(1)}_2)$-crystal:
\[
z^m (2,1,2) \otimes z^{m'} (0,2,1) \mapsto
z^{m'+3}(2,1,0) \otimes z^{m-3}(0,2,3).
\]
\end{example}

\begin{example}\label{ex:A(2)odd}
$\mathfrak{g}=A^{(2)}_7.$\\[\sukima]%

\noindent
$0: \cdots 1 1 \kakomiu{\bar{2}
\bar{4} 2 2 2} 1 1 1 1
1 1 \kakomid{\bar{2} \bar{2} \bar{2}}
1 1 1 1 1 1 1 1 1 1 1 1 1 1 1 1 1 1 1 1 1 1 1 1
1 1 1 1 1 1 1 1 1 1 1 1 1 1 \cdots$\\[\sukima]%
$1: \cdots 1 1 1 1
1 1 1 \kakomiu{\bar{2} \bar{4} 2 2
2} 1 1 1 1 \kakomid{\bar{2} \bar{2}
\bar{2}}
1 1 1 1 1 1 1 1 1 1 1 1 1 1 1 1 1 1 1 1 1
1 1 1 1 1 1 1 1 1 1 1 1 1 1 \cdots$\\[\sukima]%
$2:  \cdots 1 1 1 1
1 1 1 1 1 1 1 1 \kakomiu{\bar{2}
\bar{4} 2 2 2} 1 1 \kakomid{\bar{2}
\bar{2} \bar{2}} 1 1 1 1 1 1 1 1 1 1 1 1 1 1 1 1 1 1
1 1 1 1 1 1 1 1 1 1 1 1 1 1 \cdots$\\[\sukima]%
$3: \cdots 1 1 1 1
1 1 1 1 1 1 1 1 1
1 1 1 1 \kakomiu{\bar{2} \bar{4} 2
2 2} \kakomid{\bar{2} \bar{2} \bar{2}} 1 1 1 1 1 1 1 1 1 1 1 1 1 1 1
1 1 1 1 1 1 1 1 1 1 1 1 1 1 \cdots$\\[\sukima]%
$4: \cdots 1 1 1 1
1 1 1 1 1 1 1 1 1
1 1 1 1 1 1 1 1 1
\kakomiu{\bar{2} \bar{4} 2} \kakomiud{\bar{1} \bar{1}} \kakomid{\bar{2}}
1 1 1 1 1 1 1 1 1 1 1 1 1 1 1 1 1 1 1 1 1 1 1 1 1 1 \cdots$\\[\sukima]%
$5: \cdots 1 1 1 1
1 1 1 1 1 1 1 1 1
1 1 1 1 1 1 1 1 1
1 1 1 1 1 \kakomiu{\bar{2}} \kakomiud{\bar{1}
\bar{1} \bar{1}} \kakomiu{\bar{4}} 1 1 1 1 1 1 1 1
1 1 1 1 1 1 1 1 1 1 1 1 1 1 \cdots$\\[\sukima]%
$6: \cdots 1 1 1 1
1 1 1 1 1 1 1 1 1
1 1 1 1 1 1 1 1 1
1 1 1 1 1 1 1 1 1
\kakomid{2} \kakomiud{\bar{1} \bar{1}} \kakomiu{\bar{2} \bar{2} \bar{4}} 1 1 1
1 1 1 1 1 1 1 1 1 1 1 1 1 1 \cdots$\\[\sukima]%
$7: \cdots 1 1 1 1
1 1 1 1 1 1 1 1 1
1 1 1 1 1 1 1 1 1
1 1 1 1 1 1 1 1 1
1 1 1 \kakomid{2 2 2} \kakomiu{\bar{2} \bar{2}
\bar{2} \bar{2} \bar{4}} 1 1 1
1 1 1 1 1 1 1 1 1 \cdots$\\[\sukima]%
$8: \cdots 1 1 1 1
1 1 1 1 1 1 1 1 1
1 1 1 1 1 1 1 1 1
1 1 1 1 1 1 1 1 1
1 1 1 1 1 1 \kakomid{2 2 2}
1 1 \kakomiu{\bar{2} \bar{2} \bar{2} \bar{2} \bar{4}}
1 1 1 1 1 1 1 \cdots$\\[\sukima]%
$9: \cdots 1 1 1 1
1 1 1 1 1 1 1 1 1
1 1 1 1 1 1 1 1 1
1 1 1 1 1 1 1 1 1
1 1 1 1 1 1 1 1 1
\kakomid{2 2 2} 1 1 1 1 \kakomiu{\bar{2}
\bar{2} \bar{2} \bar{2} \bar{4}} 1 1 \cdots$\\[\sukima]%

Compare this with the combinatorial $R$  of $U_q(A^{(2)}_5)$-crystal :
\begin{align*}
&z^m (3,0,0,1,0,1) \otimes z^{m'} (0,0,0,0,0,3) \\
\mapsto
&z^{m'+0}(3,0,0,0,0,0) \otimes z^{m-0}(0,0,0,1,0,4).
\end{align*}
\end{example}

The following example starts with the formally identical initial condition
as  {\sc Example} \ref{ex:A(2)odd} but under a different $\geh$.
The outgoing solitons are different.
Moreover, separation into two final solitons
is much quicker in  {\sc Example} \ref{ex:A(2)even1}.
It has a larger phase shift, which corresponds to the greater value of
the energy functions of the combinatorial $R$.

\begin{example}\label{ex:A(2)even1}
$\mathfrak{g}=A^{(2)}_8.$\\[\sukima]%

\noindent
$0: \cdots  1 1 \kakomiu{\bar{2} \bar{4}
2 2 2} 1 1 1 1 1 1
\kakomid{\bar{2} \bar{2} \bar{2}} 1 1 1 1
1 1 1 1 1 1 1 1 1
1 1 1 1 1 1 1 1 1
1 1 1 1 1 1 1 1 1 \cdots$\\[\sukima]%
$1: \cdots  1 1 1 1 1
1 1 \kakomiu{\bar{2} \bar{4} 2 2 2}
1 1 1 1 \kakomid{\bar{2} \bar{2} \bar{2}}
1 1 1 1 1 1 1 1 1
1 1 1 1 1 1 1 1 1
1 1 1 1 1 1 1 1 1 1 \cdots$\\[\sukima]%
$2: \cdots  1 1 1 1 1
1 1 1 1 1 1 1 \kakomiu{\bar{2} \bar{4} 2 2 2} 1 1
\kakomid{\bar{2} \bar{2} \bar{2}} 1 1 1 1 1
1 1 1 1 1 1 1 1 1
1 1 1 1 1 1 1 1 1 1 1 \cdots$\\[\sukima]%
$3: \cdots  1 1 1 1 1 1 1 1 1 1 1 1 1 1
1 1 1 \kakomiu{\bar{2} \bar{4} 2 2 2} \kakomid{\bar{2} \bar{2} \bar{2}} 1 1 1
1 1 1 1 1 1 1 1 1
1 1 1 1 1 1 1 1 1 1 \cdots$\\[\sukima]%
$4: \cdots  1 1 1 1 1
1 1 1 1 1 1 1 1 1
1 1 1 1 1 1 1 1 \kakomiu{\bar{2} \bar{4} 2}
\kakomiud{\bar{1} \bar{1}} \kakomid{\bar{2}} 1
1 1 1 1 1 1 1 1 1
1 1 1 1 1 1 1 1 1 \cdots$\\[\sukima]%
$5: \cdots  1 1 1 1 1
1 1 1 1 1 1 1 1 1
1 1 1 1 1 1 1 1 1
1 1 \bar{2} \bar{4} \kakomiu{2} \kakomiud{1 1 \bar{1}}
\kakomiu{\bar{1}} \bar{2} 1 1 1 1
1 1 1 1 1 1 1 1 1 1 \cdots$\\[\sukima]%
$6: \cdots  1 1 1 1 1
1 1 1 1 1 1 1 1 1
1 1 1 1 1 1 1 1 1
1 1 1 1 1 \bar{2} \bar{4} 2
\kakomid{1} \kakomiud{1 1} \kakomiu{1 \bar{1} \bar{1}}
\bar{2} 1 1 1 1 1 1 1
1 1 \cdots$\\[\sukima]%
$7: \cdots  1 1 1 1 1
1 1 1 1 1 1 1 1 1
1 1 1 1 1 1 1 1 1
1 1 1 1 1 1 1 1 \bar{2}
\bar{4} 2 \kakomid{1 1 1} \kakomiu{1 1 1
\bar{1} \bar{1}} \bar{2} 1 1 1 1 \cdots$\\[\sukima]%

Compare this with the combinatorial $R$  of $U_q(A^{(2)}_6)$-crystal :
\begin{align*}
&z^m (3,0,0,1,0,1) \otimes z^{m'} (0,0,0,0,0,3) \\
\mapsto
&z^{m'+6}(1,0,0,1,0,1) \otimes z^{m-6}(0,0,0,0,0,1).
\end{align*}
\end{example}
Note that $A^{(2)}_8$ is Type III, so the
counting of the phase is doubled as in (\ref{eq:coef2-1}) -- (\ref{eq:coef2-2}).
In particular when $\phi \in B_1$ locates in front of solitons,
there is a special feature (\ref{eq:coef2-2}).
The following is such an example.

\begin{example}\label{ex:A(2)even2}
$\mathfrak{g}=A^{(2)}_8.$\\[\sukima]%

\noindent
$0: \cdots  1 1 \kakomiu{\bar{2} 4
3 3 3} 1 1 1 1 1 \kakomid{\emptyset
 \bar{4} 2} 1 1 1 1 1
1 1 1 1 1 1 1 1 1
1 1 1 1 1 1 1 1 1 \cdots$\\[\sukima]%
$1: \cdots  1 1 1 1 1
1 1 \kakomiu{\bar{2} 4 3 3 3} 1
1 1 \kakomid{\emptyset  \bar{4} 2} 1 1 1
1 1 1 1 1 1 1 1 1
1 1 1 1 1 1 1 1 \cdots$\\[\sukima]%
$2: \cdots  1 1 1 1 1
1 1 1 1 1 1 1 \kakomiu{\bar{2}
4 3 3 3} 1 \kakomid{\emptyset  \bar{4} 2}
1 1 1 1 1 1 1 1 1
1 1 1 1 1 1 1 1 \cdots$\\[\sukima]%
$3: \cdots  1 1 1 1 1
1 1 1 1 1 1 1 1 1
1 1 1 \kakomiu{\bar{2} 4 3 1} \kakomiud{\emptyset}
\kakomid{\bar{4} 3} 3 2 1 1 1 1
1 1 1 1 1 1 1 1 \cdots$\\[\sukima]%
$4: \cdots  1 1 1 1 1
1 1 1 1 1 1 1 1 1
1 1 1 1 1 1 \bar{2} 4
\kakomiu{3 1} \kakomiud{1 1 \emptyset}  \bar{4} 3 3
2 1 1 1 1 1 1 1 \cdots$\\[\sukima]%
$5: \cdots  1 1 1 1 1
1 1 1 1 1 1 1 1 1
1 1 1 1 1 1 1 1 1
\bar{2} 4 3 1 \kakomiud{1 1 1} \kakomiu{1
\emptyset}  \bar{4} 3 3 2 1 1 \cdots$\\[\sukima]%

Compare this with the combinatorial $R$  of $U_q(A^{(2)}_6)$-crystal :
\begin{align*}
&z^m (0,3,1,0,0,1) \otimes z^{m'} (1,0,0,1,0,0) \\
\mapsto
&z^{m'+9}(0,1,1,0,0,1) \otimes z^{m-9}(1,2,0,1,0,0).
\end{align*}
\end{example}

\begin{example}\label{ex:B(1)}
$\mathfrak{g}=B^{(1)}_4.$\\[\sukima]%

\noindent
$0: \cdots 1 1 1 1 1 1 1 \kakomiu{4 4
4 3 3} 1 1 1 1 1 1
\kakomid{\bar{4} \bar{4} 0} 1 1 1 1 1 1 1
1 1 1 1 1 1 1 1 1 1 1 1 1 1 1 1 1 1 1 1 1 1 1 1 \cdots$\\[\sukima]%
$1: \cdots 1 1 1 1 1 1 1 1 1
1 1 1 \kakomiu{4 4 4 3 3} 1
1 1 1 \kakomid{\bar{4} \bar{4} 0} 1 1 1 1 1 1 1 1 1 1 1 1 1 1
1 1 1 1 1 1 1 1 1 1 1 1 1 1 \cdots$\\[\sukima]%
$2: \cdots 1 1 1 1 1 1 1 1 1
1 1 1 1 1 1 1 1 \kakomiu{4
4 4 3 3} 1 1 \kakomid{\bar{4} \bar{4}
0} 1 1 1 1 1 1 1 1 1 1 1 1 1 1 1 1 1 1 1 1 1 1 1 1 1 \cdots$\\[\sukima]%
$3: \cdots 1 1 1 1 1 1 1 1 1
1 1 1 1 1 1 1 1 1
1 1 1 1 \kakomiu{4 4 4 3 3}
\kakomid{\bar{4} \bar{4} 0} 1 1
1 1 1 1 1 1 1 1 1 1 1 1 1 1 1 1 1 1 1 1 \cdots$\\[\sukima]%
$4: \cdots 1 1 1 1 1 1 1 1 1
1 1 1 1 1 1 1 1 1
1 1 1 1 1 1 1 1 1
\kakomiu{4 4 4} \kakomiud{\bar{4} \bar{4}} \kakomid{0} 3
3 1 1 1 1 1 1 1 1 1 1 1 1 1 1 1 1 1 \cdots$\\[\sukima]%
$5: \cdots 1 1 1 1 1 1 1 1 1
1 1 1 1 1 1 1 1 1
1 1 1 1 1 1 1 1 1
1 1 1 4 4 \kakomiu{4} \kakomiud{1 1 \bar{4}}
\kakomiu{\bar{4}} 0 3 3 1 1 1 1 1 1 1 1 1 1 1 1 \cdots$\\[\sukima]%
$6: \cdots 1 1 1 1 1 1 1 1 1
1 1 1 1 1 1 1 1 1
1 1 1 1 1 1 1 1 1
1 1 1 1 1 1 4 4 4
\kakomid{1} \kakomiud{1 1} \kakomiu{1 \bar{4} \bar{4}} 0
3 3 1 1 1 1 1 1 1 \cdots$\\[\sukima]%
$7: \cdots 1 1 1 1 1 1 1 1 1
1 1 1 1 1 1 1 1 1
1 1 1 1 1 1 1 1 1
1 1 1 1 1 1 1 1 1 4 4 4
\kakomid{1 1 1} \kakomiu{1 1 1 \bar{4} \bar{4}} 0
3 3 1 1 \cdots$\\[\sukima]%

Compare this with the combinatorial $R$  of $U_q(B^{(1)}_3)$-crystal :
\begin{align*}
&z^m (0,2,3,0,0,0,0) \otimes z^{m'} (0,0,0,1,2,0,0) \\
\mapsto
&z^{m'+3}(0,0,3,0,0,0,0) \otimes z^{m-3}(0,2,0,1,2,0,0).
\end{align*}
\end{example}

\begin{example}\label{ex:C(1)}
$\mathfrak{g}=C^{(1)}_4.$\\[\sukima]%

\noindent
$0: \cdots 1 1 \kakomiu{\bar{1}
\bar{2} \bar{3} \bar{4} 1} 1 1 1
1 1 1 \kakomid{\bar{2} \bar{4} 3} 1 1 1 1 1 1 1 1 1 1 1 1 1 1
 1 1 1 1 1 1 1 1 1 1 1 1 1 \cdots$\\[\sukima]%
$1: \cdots 1 1 1 1
1 1 1 \kakomiu{\bar{1} \bar{2} \bar{3} \bar{4}
1} 1 1 1 1 \kakomid{\bar{2} \bar{4}
3} 1 1 1 1 1 1 1 1 1 1 1 1 1 1 1 1 1 1 1 1 1 1 1 1 \cdots$\\[\sukima]%
$2: \cdots 1 1 1 1
1 1 1 1 1 1 1 1 \kakomiu{\bar{1}
\bar{2} \bar{3} \bar{4} 1} 1 1 \kakomid{\bar{2}
\bar{4} 3} 1 1 1 1 1 1 1 1 1 1 1 1 1 1 1 1 1 1 1 1 1 \cdots$\\[\sukima]%
$3: \cdots 1 1 1 1
1 1 1 1 1 1 1 1 1
1 1 1 1 \kakomiu{\bar{1} \bar{2} \bar{3}
1 1} \kakomid{\bar{2} \bar{4} \bar{4}} 3
 1 1 1 1 1 1 1 1 1 1 1 1 1 1 1 1 1 \cdots$\\[\sukima]%
$4: \cdots 1 1 1 1
1 1 1 1 1 1 1 1 1
1 1 1 1 1 1 1 1 \bar{1}
\kakomiu{\bar{3} 1 1} \kakomiud{1 \bar{2}} \kakomid{\bar{2}} \bar{4}
\bar{4} 3 1 1 1 1 1 1 1 1 1 1 1 1 \cdots$\\[\sukima]%
$5: \cdots 1 1 1 1
1 1 1 1 1 1 1 1 1
1 1 1 1 1 1 1 1 1
1 1 \bar{1} \bar{3} 1 \kakomiu{1} \kakomiud{1
1 1} \kakomiu{\bar{2}} \bar{2} \bar{4} \bar{4} 3 1 1 1 1 1 1 1 \cdots$\\[\sukima]%
$6: \cdots 1 1 1 1
1 1 1 1 1 1 1 1 1
1 1 1 1 1 1 1 1 1
1 1 1 1 1 \bar{1} \bar{3} 1 1 \kakomid{1}
\kakomiud{1 1} \kakomiu{1 1 \bar{2}} \bar{2} \bar{4} \bar{4} 3 1 1 \cdots$\\[\sukima]%

Compare this with the combinatorial $R$  of $U_q(C^{(1)}_3)$-crystal :
\begin{align*}
&z^m (0,0,0,1,1,1) \otimes z^{m'} (0,1,0,1,0,1) \\ \mapsto
&z^{m'+4}(0,0,0,0,1,0) \otimes z^{m-4}(0,1,0,2,0,2).
\end{align*}
\end{example}

\begin{example}\label{ex:D(1)}
$\mathfrak{g}=D^{(1)}_5.$\\[\sukima]%

\noindent
$0: \cdots 1 1 \kakomiu{\bar{3} \bar{5} 3 2 2} 1 1 1 1 1 1
\kakomid{\bar{2} \bar{2} \bar{5}} 1 1 1 1 1 1 1 1 1 1 1 1 1 1 1 1
1 1 1 1 1 1 1 1 1 1 1 1 1 1 1 1 1 1 1 1 1 1 1 \cdots$\\[\sukima]%
$1: \cdots 1 1 1 1 1 1 1 \kakomiu{\bar{3} \bar{5} 3 2 2} 1 1 1 1
\kakomid{\bar{2} \bar{2} \bar{5}} 1 1 1 1 1 1 1 1 1 1 1 1 1
1 1 1 1 1 1 1 1 1 1 1 1 1 1 1 1 1 1 1 1 1 1 1 \cdots$\\[\sukima]%
$2: \cdots 1 1 1 1 1 1 1 1 1 1 1 1
\kakomiu{\bar{3} \bar{5} 3 2 2} 1 1 \kakomid{\bar{2} \bar{2} \bar{5}}
 1 1 1 1 1 1 1 1 1 1 1 1 1 1 1 1 1 1 1 1 1 1 1 1 1 1 1 1 1 1 1 1 1 \cdots$\\[\sukima]%
$3: \cdots 1 1 1 1 1 1 1 1 1 1 1 1 1 1 1 1 1
\kakomiu{\bar{3} \bar{5} 3 2 2} \kakomid{\bar{2} \bar{2} \bar{5}}
 1 1 1 1 1 1 1 1 1 1 1 1 1 1 1 1 1 1 1 1 1 1 1 1 1 1 1 1 1 1 \cdots$\\[\sukima]%
$4: \cdots 1 1 1 1 1 1 1 1 1 1 1 1 1 1 1 1 1 1 1 1
 1 1 \kakomiu{\bar{3} \bar{5} 3} \kakomiud{\bar{1} \bar{1}}
 \kakomid{\bar{5}} 1 1 1 1 1 1 1 1 1 1 1 1 1 1 1 1 1 1 1 1
 1 1 1 1 1 1 1 \cdots$\\[\sukima]%
$5: \cdots 1 1 1 1 1 1 1 1 1 1 1 1 1 1 1 1 1 1 1 1
 1 1 1 1 1 1 1 \kakomiu{4} \kakomiud{\bar{1} \bar{1} \bar{4}}
 \kakomiu{\bar{5}} \bar{5}
 1 1 1 1 1 1 1 1 1 1 1 1 1 1 1 1 1 1 1 1 1 1 \cdots$\\[\sukima]%
$6: \cdots 1 1 1 1 1 1 1 1 1 1 1 1 1 1 1 1 1 1 1 1
 1 1 1 1 1 1 1 1 1 1 4 \kakomid{2} \kakomiud{2
\bar{2}} \kakomiu{\bar{2} \bar{4} \bar{5}} \bar{5}
1 1 1 1 1 1 1 1 1 1 1 1 1 1 1 1 1 \cdots$\\[\sukima]%
$7: \cdots 1 1 1 1 1 1 1 1 1 1 1 1 1 1 1 1 1 1 1 1
 1 1 1 1 1 1 1 1 1 1 1 1 1 4 \kakomid{2 2 1} \kakomiu{1
\bar{2} \bar{2} \bar{4} \bar{5}} \bar{5} 1 1 1 1 1 1 1 1 1 1 1 1 \cdots$\\[\sukima]%
$8: \cdots 1 1 1 1 1 1 1 1 1 1 1 1 1 1 1 1 1 1 1 1
 1 1 1 1 1 1 1 1 1 1 1 1 1 1 1 1 4 \kakomid{2 2 1} 1 1 \kakomiu{1
\bar{2} \bar{2} \bar{4} \bar{5}} \bar{5} 1 1 1 1 1 1 1 \cdots$\\[\sukima]%
$9: \cdots 1 1 1 1 1 1 1 1 1 1 1 1 1 1 1 1 1 1 1 1
 1 1 1 1 1 1 1 1 1 1 1 1 1 1 1 1 1 1 1 4 \kakomiu{2 2 1} 1 1 1 1
 \kakomiu{1 \bar{2} \bar{2} \bar{4} \bar{5}} \bar{5} 1 1 \cdots$\\[\sukima]%

Compare this with the combinatorial $R$  of $U_q(D^{(1)}_4)$-crystal :
\begin{align*}
&z^m (2,1,0,0,1,0,1,0) \otimes z^{m'} (0,0,0,0,1,0,0,2) \\ \mapsto
&z^{m'+1}(2,0,1,0,0,0,0,0) \otimes z^{m-1}(0,0,0,0,2,1,0,2).
\end{align*}
\end{example}

\begin{example}\label{ex:D(2)}
$\mathfrak{g}=D^{(2)}_5.$\\[\sukima]%

\noindent
$0: \cdots 1 1 \kakomiu{\bar{2}
0 3 2 2} 1 1 1 1 1
1 \kakomid{\bar{1} 0 1} 1 1 1 1
 1 1 1 1 1 1 1 1 1 1 1 1 1 1 1 1 1 1 1 1 1 1 1 \cdots$\\[\sukima]%
$1: \cdots 1 1 1 1
1 1 1 \kakomiu{\bar{2} 0 3 2 2}
1 1 1 1 \kakomid{\bar{1} 0 1}
1 1 1 1 1 1 1 1 1 1 1 1 1 1 1 1 1 1 1 1 1 1 1 1 \cdots$\\[\sukima]%
$2: \cdots 1 1 1 1
1 1 1 1 1 1 1 1 \kakomiu{\bar{2}
0 3 2 2} 1 1 \kakomid{\bar{1} 0 1} 1
1 1 1 1 1 1 1 1 1 1 1 1 1 1 1 1 1 1 1 1  \cdots$\\[\sukima]%
$3: \cdots 1 1 1 1
1 1 1 1 1 1 1 1 1
1 1 1 1 \kakomiu{\bar{2} 0 3 2
1} \kakomid{\bar{1} 0 2} 1 1 1 1 1 1 1 1 1 1 1 1 1 1 1 1 1 1 \cdots$\\[\sukima]%
$4: \cdots 1 1 1 1
1 1 1 1 1 1 1 1 1
1 1 1 1 1 1 1 1 \bar{2}
\kakomiu{0 3 1} \kakomiud{1 \bar{1}} \kakomid{0}
2 2 1 1 1 1 1 1 1 1 1 1 1 1 1 \cdots$\\[\sukima]%
$5: \cdots 1 1 1 1
1 1 1 1 1 1 1 1 1
1 1 1 1 1 1 1 1 1
1 1 \bar{2} 0 3 \kakomiu{1} \kakomiud{1 1
1} \kakomiu{\bar{1}} 0 2 2 1 1 1 1 1 1 1 1 \cdots$\\[\sukima]%
$6: \cdots 1 1 1 1
1 1 1 1 1 1 1 1 1
1 1 1 1 1 1 1 1 1
1 1 1 1 1 \bar{2} 0 3 1 \kakomid{1} \kakomiud{1 1} \kakomiu{1
1 \bar{1}} 0 2 2 1 1 1 \cdots$\\[\sukima]%

Compare this with the combinatorial
$R$  of $U_q(D^{(2)}_4)$-crystal :
\begin{align*}
&z^m (2,1,0,1,0,0,1) \otimes z^{m'} (0,0,0,1,0,0,0) \\ \mapsto
&z^{m'+8}(0,1,0,1,0,0,1) \otimes z^{m-8}(2,0,0,1,0,0,0).
\end{align*}
\end{example}
Note that $D^{(2)}_5$ is Type III, 
so the counting of the phase is doubled as in 
{\sc Examples} \ref{ex:A(2)even1} and \ref{ex:A(2)even2}.

Finally we present two comparable figures showing the
soliton nature, namely, independence of the order of the collisions.

\begin{example}
$\mathfrak{g}=B^{(1)}_4.$\\[\sukima]%

\noindent
$\makebox[1em][r]{$0$}: \cdots  1 1 \kakomiu{\bar{2} \bar{2}
2} 1 1 1 \kakomid{\bar{3} \bar{3}} 1
1 1 1 1 1 1 1 \kakomis{0} 1
1 1 1 1 1 1 1 1 1
1 1 1 1 1 1 1 1 1
1 1 1 1 1 1 1 1 1
1 \cdots$\\[\sukima]%
$\makebox[1em][r]{$1$}: \cdots  1 1 1 1 1
\kakomiu{\bar{2} \bar{2} 2} 1 1 \kakomid{\bar{3} \bar{3}}
1 1 1 1 1 1 1 \kakomis{0} 1
1 1 1 1 1 1 1 1 1
1 1 1 1 1 1 1 1 1
1 1 1 1 1 1 1 1 1
\cdots$\\[\sukima]%
$\makebox[1em][r]{$2$}: \cdots  1 1 1 1 1
1 1 1 \kakomiu{\bar{2} \bar{2} 2} 1
\kakomid{\bar{3} \bar{3}} 1 1 1 1 1
1 \kakomis{0} 1 1 1 1 1 1 1
1 1 1 1 1 1 1 1 1
1 1 1 1 1 1 1 1 1
1 1 \cdots$\\[\sukima]%
$\makebox[1em][r]{$3$}: \cdots  1 1 1 1 1
1 1 1 1 1 1 \kakomiu{\bar{2} 3
1} \kakomid{\bar{3} \bar{3}} \bar{3} 1 1 1
1 \kakomis{0} 1 1 1 1 1 1 1
1 1 1 1 1 1 1 1 1
1 1 1 1 1 1 1 1 1
1 \cdots$\\[\sukima]%
$\makebox[1em][r]{$4$}: \cdots  1 1 1 1 1
1 1 1 1 1 1 1 1 \bar{2}
\kakomiu{3 1} \kakomiud{1} \kakomid{\bar{3}} \bar{3} \bar{3} 1
1 \kakomis{0} 1 1 1 1 1 1 1
1 1 1 1 1 1 1 1 1
1 1 1 1 1 1 1 1 1
\cdots$\\[\sukima]%
$\makebox[1em][r]{$5$}: \cdots  1 1 1 1 1
1 1 1 1 1 1 1 1 1
1 \bar{2} 3 \kakomiu{1} \kakomiud{1 1} \bar{3}
\bar{3} 1 \kakomis{\bar{3}} 0 1 1 1
1 1 1 1 1 1 1 1 1
1 1 1 1 1 1 1 1 1
1 1 \cdots$\\[\sukima]%
$\makebox[1em][r]{$6$}: \cdots  1 1 1 1 1
1 1 1 1 1 1 1 1 1
1 1 1 \bar{2} 3 1 \kakomiud{1 1}
\kakomid{\bar{3}} 1 \kakomis{1} \bar{3} \bar{3} 0 1
1 1 1 1 1 1 1 1 1
1 1 1 1 1 1 1 1 1
1 \cdots$\\[\sukima]%
$\makebox[1em][r]{$7$}: \cdots  1 1 1 1 1
1 1 1 1 1 1 1 1 1
1 1 1 1 1 \bar{2} 3 1
\kakomid{1} \kakomiud{\bar{3}} \kakomiu{1} \kakomius{1} 1 1 \bar{3}
\bar{3} 0 1 1 1 1 1 1
1 1 1 1 1 1 1 1 1
1 1 \cdots$\\[\sukima]%
$\makebox[1em][r]{$8$}: \cdots  1 1 1 1 1
1 1 1 1 1 1 1 1 1
1 1 1 1 1 1 1 \bar{2}
3 1 \kakomid{\bar{3} 1} \kakomius{1} \kakomiu{1 1} 1
1 \bar{3} \bar{3} 0 1 1 1
1 1 1 1 1 1 1 1 1
1 1 \cdots$\\[\sukima]%
$\makebox[1em][r]{$9$}: \cdots  1 1 1 1 1
1 1 1 1 1 1 1 1 1
1 1 1 1 1 1 1 1 1
\bar{2} 3 \bar{3} \kakomid{1} \kakomids{1} 1 \kakomiu{1
1 1} 1 1 \bar{3} \bar{3} 0
1 1 1 1 1 1 1 1 1
1 1 \cdots$\\[\sukima]%
$\makebox[1em][r]{$10$}: \cdots  1 1 1 1 1
1 1 1 1 1 1 1 1 1
1 1 1 1 1 1 1 1 1
1 1 3 \bar{2} \bar{3} \kakomids{1} \kakomid{1}
1 1 \kakomiu{1 1 1} 1 1 \bar{3}
\bar{3} 0 1 1 1 1 1 1
1 1 \cdots$\\[\sukima]%
$\makebox[1em][r]{$11$}: \cdots  1 1 1 1 1
1 1 1 1 1 1 1 1 1
1 1 1 1 1 1 1 1 1
1 1 1 3 1 \bar{2} \kakomis{\bar{3}}
\kakomid{1 1} 1 1 1 \kakomiu{1 1 1} 1
1 \bar{3} \bar{3} 0 1 1 1
1 1 \cdots$\\[\sukima]%
$\makebox[1em][r]{$12$}: \cdots  1 1 1 1 1
1 1 1 1 1 1 1 1 1
1 1 1 1 1 1 1 1 1
1 1 1 1 3 1 1 \kakomis{\bar{2}}
\bar{3} \kakomid{1 1} 1 1 1 1 \kakomiu{1
1 1} 1 1 \bar{3} \bar{3} 0
1 1 \cdots$\\[\sukima]%

\noindent
$\makebox[1em][r]{$0$}: \cdots  1 1 \kakomiu{\bar{2} \bar{2}
2} 1 1 1 1 1 1 1 1
1 1 \kakomid{\bar{3} \bar{3}} 1 1 1
\kakomis{0} 1 1 1 1 1 1 1 1
1 1 1 1 1 1 1 1 1
1 1 1 1 1 1 1 1 1
1 \cdots$\\[\sukima]%
$\makebox[1em][r]{$1$}: \cdots  1 1 1 1 1
\kakomiu{\bar{2} \bar{2} 2} 1 1 1 1
1 1 1 1 1 \kakomid{\bar{3} \bar{3}}
1 1 \kakomis{0} 1 1 1 1 1 1
1 1 1 1 1 1 1 1 1
1 1 1 1 1 1 1 1 1
1 1 \cdots$\\[\sukima]%
$\makebox[1em][r]{$2$}: \cdots  1 1 1 1 1
1 1 1 \kakomiu{\bar{2} \bar{2} 2} 1
1 1 1 1 1 1 1 \kakomid{\bar{3} \bar{3}} 1
\kakomis{0} 1 1 1 1 1
1 1 1 1 1 1 1 1 1
1 1 1 1 1 1 1 1 1
1 1 \cdots$\\[\sukima]%
$\makebox[1em][r]{$3$}: \cdots  1 1 1 1 1
1 1 1 1 1 1 \kakomiu{\bar{2} \bar{2}
2} 1 1 1 1 1 1 1 \kakomid{\bar{3}
1} \kakomis{\bar{3}} 0 1 1 1 1 1
1 1 1 1 1 1 1 1 1
1 1 1 1 1 1 1 1 1
\cdots$\\[\sukima]%
$\makebox[1em][r]{$4$}: \cdots  1 1 1 1 1
1 1 1 1 1 1 1 1 1
\kakomiu{\bar{2} \bar{2} 2} 1 1 1 1
1 \bar{3} \kakomid{1} \kakomids{1} \bar{3} 0 1
1 1 1 1 1 1 1 1 1
1 1 1 1 1 1 1 1 1
1 1 \cdots$\\[\sukima]%
$\makebox[1em][r]{$5$}: \cdots  1 1 1 1 1
1 1 1 1 1 1 1 1 1
1 1 1 \kakomiu{\bar{2} \bar{2} 2} 1
1 1 \bar{3} 1 \kakomids{1} \kakomid{1} \bar{3}
0 1 1 1 1 1 1 1 1
1 1 1 1 1 1 1 1 1
1 1 \cdots$\\[\sukima]%
$\makebox[1em][r]{$6$}: \cdots  1 1 1 1 1
1 1 1 1 1 1 1 1 1
1 1 1 1 1 1 \kakomiu{\bar{2} \bar{2}
2} 1 \bar{3} 1 \kakomis{1} \kakomid{1 1} \bar{3}
0 1 1 1 1 1 1 1 1
1 1 1 1 1 1 1 1 1
\cdots$\\[\sukima]%
$\makebox[1em][r]{$7$}: \cdots  1 1 1 1 1
1 1 1 1 1 1 1 1 1
1 1 1 1 1 1 1 1 1
\kakomiu{\bar{2} 3 \bar{3}} \bar{3} \kakomis{1} 1 \kakomid{1
1} \bar{3} 0 1 1 1 1 1
1 1 1 1 1 1 1 1 1
1 \cdots$\\[\sukima]%
$\makebox[1em][r]{$8$}: \cdots  1 1 1 1 1
1 1 1 1 1 1 1 1 1
1 1 1 1 1 1 1 1 1
1 1 3 \kakomiu{1 \bar{2}} \kakomius{\bar{3}} \bar{3}
1 \kakomid{1 1} \bar{3} 0 1 1 1
1 1 1 1 1 1 1 1 1
1 \cdots$\\[\sukima]%
$\makebox[1em][r]{$9$}: \cdots  1 1 1 1 1
1 1 1 1 1 1 1 1 1
1 1 1 1 1 1 1 1 1
1 1 1 3 1 1 \kakomius{1} \kakomiu{\bar{2}
\bar{3}} \bar{3} \kakomid{1 1} \bar{3} 0 1
1 1 1 1 1 1 1 1 1
1 \cdots$\\[\sukima]%
$\makebox[1em][r]{$10$}: \cdots  1 1 1 1 1
1 1 1 1 1 1 1 1 1
1 1 1 1 1 1 1 1 1
1 1 1 1 3 1 1 \kakomis{1} 1
\kakomiu{1 \bar{2} \bar{3}} \kakomid{1 1} \bar{3} \bar{3}
0 1 1 1 1 1 1 1 1
\cdots$\\[\sukima]%
$\makebox[1em][r]{$11$}: \cdots  1 1 1 1 1
1 1 1 1 1 1 1 1 1
1 1 1 1 1 1 1 1 1
1 1 1 1 1 3 1 1 \kakomis{1}
1 1 1 \kakomiu{\bar{2} \bar{3}} \kakomiud{1} \kakomid{1}
1 \bar{3} \bar{3} 0 1 1 1
1 1 \cdots$\\[\sukima]%
$\makebox[1em][r]{$12$}: \cdots  1 1 1 1 1
1 1 1 1 1 1 1 1 1
1 1 1 1 1 1 1 1 1
1 1 1 1 1 1 3 1 1
\kakomis{1} 1 1 1 1 \bar{2} \kakomiu{\bar{3}}
\kakomiud{1 1} 1 1 \bar{3} \bar{3} 0
1 1 \cdots$\\[\sukima]%

The both  scatterings are described by
\begin{align*}
&z^m (1,0,0,0,0,0,2) \otimes z^{m'} (0,0,0,0,0,2,0) \otimes z^{m''} (0,0,0,1,0,0,0) \\
&\mapsto
z^{m''+3} (0,1,0,0,0,0,0) \otimes z^{m'+2}(0,0,0,0,0,1,1)
\otimes z^{m-5}(0,0,0,1,0,2,0)
\end{align*}
under the combinatorial $R$  of $U_q(B^{(1)}_3)$-crystal :
$\Aff(B_3) \otimes \Aff(B_2) \otimes \Aff(B_1) \xrightarrow{\sim}
\Aff(B_1) \otimes \Aff(B_2) \otimes \Aff(B_3)$.

\end{example}

\section{Discussion}\label{sec:discussion}

So far we have only dealt with the homogeneous tensor product
$\cd \ot B_1 \ot B_1\ot \cd$ with a ferromagnetic boundary condition 
as the automaton states.
It is natural to consider a generalization such that 
the states $\cd \ot b_j \ot b_{j+1} \ot \cd$ 
are taken {}from the inhomogeneous 
tensor product $\cd \ot B_{m_j} \ot B_{m_{j+1}}\ot \cd$,
where $m_j$'s are arbitrary positive integers.
See \cite{HHIKTT} for $A^{(1)}_n$ case.
Imposing the boundary condition $b_j = u_{m_j}$ (see (\ref{eq:udef}))
for all the distant $j$, the time evolutions can be defined 
in an analogous manner to Section \ref{subsection:states}.
Let us discuss solitons and their scattering rule in 
the inhomogeneous case.

First, in order to define solitons, we actually consider 
the tensor product of the form:
\[
\cd \ot B_1  \ot B_1 \ot B_{m_1} \ot \cd \ot B_{m_N}
\ot B_1 \ot  B_1 \ot \cd.
\]
Here the left region $\cd \ot B_1 \ot B_1$ and the right region 
$B_1 \ot  B_1 \ot \cd$ are supplemented as the detectors of 
incoming and outgoing solitons.
In those regions one can assign elements in 
$\Aff(B_{l_1}) \ot \cd \ot \Aff(B_{l_i})$
to any asymptotic (enough separated) $i$-soliton states
in the same way as (\ref{eq:associate}).
One caution necessary here is how to specify the phase of 
solitons in the left region.
To count the ``distance"  $\gamma(t)$ {}from the 
right (see after (\ref{eq:sol})),
we understand that each component $B_{m_j}$ in the middle
region is effectively occupying $m_j$ slots.
In this way, to a soliton of length $l$ in either asymptotic region,
one can assign the data $z^\gamma b$ with $b$ {}from the 
$U'_q(\geh_{n-1})$-crystal $B_l$.

Now we are ready to study the collisions of solitons 
which take place in the 
middle inhomogeneous region in general.
Consider the incoming two-soliton state represented by 
\[
z^{\gamma_1}b_1 \ot z^{\gamma_2}b_2 \in 
\Aff(B_l) \ot \Aff(B_k), \quad (l > k)
\]
in the above sense, and let 
$z^{\gamma_2+\delta}{\tilde b}_2 \ot 
z^{\gamma_1-\delta}{\tilde b}_1 \in 
\Aff(B_k) \ot \Aff(B_l)$ with $\delta = H(b_1 \ot b_2)$ 
be the image under the combinatorial 
$R: \Aff(B_l) \ot \Aff(B_k) \stackrel{\sim}{\rightarrow}
\Aff(B_k) \ot \Aff(B_l)$.
Then we claim that after sufficiently many time steps under 
the evolution $T_r (r > k)$, the state again becomes the 
two-soliton state corresponding to 
\begin{align*}
&z^{\gamma_2+\delta-\Delta_k}{\tilde b}_2 \ot 
z^{\gamma_1-\delta-\Delta_l}{\tilde b}_1 \in 
\Aff(B_k) \ot \Aff(B_l), \\
&\Delta_i = \sum_{j=1}^N(m_j-i)_+.
\end{align*}
Comparing this with {\sc Theorem} \ref{th:main},
we see that the effect of inhomogeneity 
only enters the additional phase shifts 
$\Delta_l$ and $\Delta_k$.
In fact the claim reduces to {\sc Theorem} \ref{th:main} 
when $\forall m_j = 1$.
To justify the claim let us 
first note that for any $b \in B_i$ and $L, L' \in \Z_{\ge 0}$,
the relation 
\begin{align*}
&\hako{1}^{\ot(L+X_i)} \ot \imath_i(b)\ot \hako{1}^{\ot L'} \ot 
u_{m_1} \ot \cd \ot u_{m_N}\\
&\simeq 
u_{m_1} \ot \cd \ot u_{m_N} \ot \hako{1}^{\ot L} \ot \imath_i(b) 
\ot \hako{1}^{\ot (L'+X_i)}
\end{align*}
is valid, where 
$X_i = \sum_{j=1}^N\min(i,m_j) (= \sum_{j=1}^N m_j - \Delta_i)$.
This can be shown by considering the evolution of the 
soliton $\imath_i(b)$ under $T_{m_1}\cdots T_{m_N}$ and applying
{\sc Lemma} \ref{lem:single-soliton} (2).
It means that the soliton $\imath_i(b)$ gains the phase shift 
$L'+X_i - (L' + \sum_{j=1}^Nm_j) = -\Delta_i$ when penetrating the 
inhomogeneous region $u_{m_1} \ot \cd \ot u_{m_N}$.
We depict the result as
\begin{align}
&\hskip8.7cm {\scriptstyle -\Delta_i} \nonumber\\
&....\imath_i(b)....u_{m_1}\cdots u_{m_N}............\; 
\stackrel{T_r^t}{\rightarrow}\; 
.............u_{m_1}\cdots u_{m_N}....\imath_i(b).... \label{eq:ps}
\end{align}
for some $t \gg 1$, where $...$ denotes the tensor products of 
$\hako{1} \in B_1$.

Now we proceed to the consideration on 
the incoming state in the claim.
In the above notation it looks as
\[
p_1 = 
.....{\imath}_l(b_1).....{\imath}_k(b_2).........
u_{m_1}\cdots u_{m_N}.................................................. \;,
\]
where we have suppressed the initial phase $\gamma_1, \gamma_2$.
We shall keep track of their
shifts only.
The collision may take place in the middle region.
However by means of (\ref{eq:uu}) 
one can push  it away to the right as far as one wishes.
Namely $p_1 \simeq p_1'$, where
\[
p_1' = 
.....{\imath}_l(b_1).....{\imath}_k(b_2).............
...............u_{m_1}\cdots u_{m_N}...
............................ \; .
\]
Pushing it away enough, one can let two solitons 
finish the collision within the enlarged left region.
To such a collision {\sc Theorem} \ref{th:main} is applied, leading to 
\begin{align*}
&\hskip4.1cm {\scriptstyle \delta}\hskip1.2cm {\scriptstyle -\delta} \\
&T_r^{t_1}(p'_1) = 
.......................
{\imath}_k(\tilde{b}_2).....{\imath}_l(\tilde{b}_1)..
..u_{m_1}\cdots u_{m_N}..........
........................
\end{align*}
for some  $t_1 \gg 1$, where 
the phase shift is given by $\delta = H(b_1 \ot b_2)$.
Now we can utilize (\ref{eq:ps}) to let the two solitons 
penetrate  the inhomogeneous region. 
This can be done individually for each soliton because the 
two solitons in the 
above may be assumed to be enough separated and the right one 
${\imath}_l(\tilde{b}_1)$ is faster than 
${\imath}_k(\tilde{b}_2)$ under $T_r$ with $r>k$.
Thus for some $t_2 \gg t_1$ we have 
\begin{align*}
&\hskip9.1cm {\scriptstyle \delta-\Delta_k}
\hskip0.4cm {\scriptstyle -\delta-\Delta_l} \\
&T_r^{t_2}(p'_1) = 
.............................................
.....u_{m_1}\cdots u_{m_N}.......
{\imath}_k(\tilde{b}_2).....{\imath}_l(\tilde{b}_1).... \; ,
\end{align*}
where the solitons have gained the extra phase shifts
$-\Delta_l$ and $-\Delta_k$.
Finally we can pull the region $u_{m_1}\cdots u_{m_N}$ back to
the original position again by (\ref{eq:uu}).
As the result we have $T_r^{t_2}(p'_1) \simeq p_2$, where 
\begin{align*}
&\hskip9.0cm {\scriptstyle \delta-\Delta_k}
\hskip0.4cm {\scriptstyle -\delta-\Delta_l} \\
&p_2 = 
.................................
.....u_{m_1}\cdots u_{m_N}..........................
{\imath}_k(\tilde{b}_2).....{\imath}_l(\tilde{b}_1).... \; ,
\end{align*}
and this tensor product has the same ordering  of the components 
$B_1$'s and $B_{m_j}$'s  as $p_1$.
Since $T_r$ commutes with the isomorphism $\simeq$,
we conclude $T_r^{t_2}(p_1) = p_2$.
This completes the proof of the claim.

\begin{example}\label{ex:inhomo}
Consider the  $\mathfrak{g}=C^{(1)}_3$ automaton 
with the inhomogeneous region $m_1 = 3, m_2 = 2 \, (N=2)$.
The following is the collision of two solitons with length $l=3$ and $k=1$
in the inhomogeneous region under the time evolution $T_{12}$:

\vskip0.3cm
\noindent
$\kakomiu{\bar{2} \cdot  \bar{3} \cdot  3} \cdot
 1 \cdot  1 \cdot  1 \cdot  1
\cdot  \kakomid{2} \cdot  1 1 1 \cdot  1
1 \cdot  1 \cdot  1 \cdot  1 \cdot
 1 \cdot  1 \cdot  1 \cdot  1
\cdot  1 \cdot  1 \cdot  1 \cdot 
1 \cdot  1$\vskip0.1cm\noindent
$1 \cdot  1 \cdot  1 \cdot  \kakomiu{\bar{2}
\cdot  \bar{3} \cdot  3} \cdot  1 \cdot
 1 \cdot  \kakomid{1} 1 2 \cdot  1
1 \cdot  1 \cdot  1 \cdot  1 \cdot
 1 \cdot  1 \cdot  1 \cdot  1
\cdot  1 \cdot  1 \cdot  1 \cdot 
1 \cdot  1$\vskip0.1cm\noindent
$1 \cdot  1 \cdot  1 \cdot  1
\cdot  1 \cdot  1 \cdot  \kakomiu{\bar{2} \cdot
 \bar{3} \cdot  1} \kakomid{1} 3 \cdot 
1 2 \cdot  1 \cdot  1 \cdot  1
\cdot  1 \cdot  1 \cdot  1 \cdot 
1 \cdot  1 \cdot  1 \cdot  1 \cdot
 1 \cdot  1$\vskip0.1cm\noindent
$1 \cdot  1 \cdot  1 \cdot  1
\cdot  1 \cdot  1 \cdot  1 \cdot 
1 \cdot  1 \kakomiu{\bar{3}} \kakomiud{\bar{2}} \kakomiu { \; \cdot  \, 1}
3 \cdot  2 \cdot  1 \cdot  1 \cdot
 1 \cdot  1 \cdot  1 \cdot  1
\cdot  1 \cdot  1 \cdot  1 \cdot 
1 \cdot  1$\vskip0.1cm\noindent
$1 \cdot  1 \cdot  1 \cdot  1
\cdot  1 \cdot  1 \cdot  1 \cdot 
1 \cdot  1 1 1  \cdot  \kakomid{1} \kakomiu{\bar{2}
\cdot  1 \cdot  \bar{3}} \cdot  3 \cdot
 2 \cdot  1 \cdot  1 \cdot  1
\cdot  1 \cdot  1 \cdot  1 \cdot 
1 \cdot  1$\vskip0.1cm\noindent
$1 \cdot  1 \cdot  1 \cdot  1
\cdot  1 \cdot  1 \cdot  1 \cdot 
1 \cdot  1 1 1  \cdot  1 \kakomid{1}
\cdot  \bar{2} \cdot  1 \cdot  \kakomiu{1 \cdot
 1 \cdot  \bar{3}} \cdot  3 \cdot 
2 \cdot  1 \cdot  1 \cdot  1 \cdot
 1 \cdot  1$\vskip0.1cm\noindent
$1 \cdot  1 \cdot  1 \cdot  1
\cdot  1 \cdot  1 \cdot  1 \cdot 
1 \cdot  1 1 1 \cdot  1 1
\cdot  \kakomid{1} \cdot  \bar{2} \cdot  1 \cdot
 1 \cdot  1 \cdot  \kakomiu{1 \cdot  1
\cdot  \bar{3}} \cdot  3 \cdot  2 \cdot
 1 \cdot  1$\\[1em]
\noindent
Here $\cdot$ stands for $\otimes$, and 
$111:=(3,0,0,0,0,0)\in B_3,\,
112:=(2,1,0,0,0,0)\in B_3,\,
1 \bar{3} \bar{2}:=(1,0,0,1,1,0) \in B_3,
1 \bar{2} :=(1,0,0,0,1,0) \in B_2$, etc.
On the other hand, one has 
\[
z^m (0,1,1,1) \otimes z^{m'} (1,0,0,0) \mapsto
z^{m'+2}(0,0,0,1) \otimes z^{m-2}(1,1,1,0)
\]
under the combinatorial $R$  of $U_q(C^{(1)}_2)$-crystal :
$\Aff(B_3) \otimes \Aff(B_1) \xrightarrow{\sim}
\Aff(B_1) \otimes \Aff(B_3)$.
Since $\Delta_3=0$ and $\Delta_1=3$, the phase shifts of the 
larger and smaller solitons should be 
$-2-\Delta_3=-2$ and $2-\Delta_1=-1$, respectively.
These values indeed agree with the deviation of the outgoing solitons 
{}from the corresponding markers.
\end{example}

\vspace{0.4cm}
\noindent
{\bf Acknowledgements}. \hspace{0.1cm}
The authors thank M. Kashiwara for useful discussion.
A.K. and M.O. thank 
S. Berman, P. Fendley, Y. Huang, K. Misra, and B. Parshall 
for the warm hospitality at 
{\it Infinite-Dimensional Lie Theory and Conformal Field Theory}, 
held at University of Virginia,
Charlottesville during  May 23--27, 2000, where a part of this work 
was presented.

\appendix

\section{\mathversion{bold} Crystals and the combinatorial $R$ }\label{section:CCR}

We review fundamental facts on crystals. It is assumed that our crystals
are crystal bases of irreducible finite-dimensional representations of a 
quantum affine algebra $U'_q(\geh)$. Our basic reference is \cite{KMN}.

Let $P$ be the weight lattice, $\{\alpha_i\}_{0\le i\le n}$ the simple 
roots, and $\{\La_i\}_{0\le i\le n}$ the fundamental weights of $\geh$.
Let $\Lab_i$ denote the classical part of $\La_i$. See \cite{K} for details.
A crystal $B$ is a finite set with weight decomposition $B=\sqcup_{\la\in P}
B_\la$. The Kashiwara operators $\et{i},\ft{i}$ ($i=0,1,\cd,n$) act on
$B$ as 
\[
\et{i}: B_\la\longrightarrow B_{\la+\alpha_i} \sqcup \{0\},\quad
\ft{i}: B_\la\longrightarrow B_{\la-\alpha_i} \sqcup \{0\}.
\]
In particular, these operators are nilpotent.
By definition, we have $\ft{i}b=b'$ if and only if $b=\et{i}b'$.
Drawing $b\stackrel{i}{\rightarrow}b'$ in such a case, $B$ is endowed 
with the structure of colored oriented graph called {\em crystal graph}.
For any $b\in B$, set 
\[
\veps_i(b)=\max\{m\ge0\mid \et{i}^m b\ne0\},\quad
\vphi_i(b)=\max\{m\ge0\mid \ft{i}^m b\ne0\}.
\]
Then we have the weight $\wt b$ of $b$ by $\wt b=\sum_{i=0}^n(\vphi_i(b)-
\veps_i(b))\La_i$.

For two crystals $B$ and $B'$, one can define the tensor product
$B\ot B'=\{b\ot b'\mid b\in B,b'\in B'\}$. The operators $\et{i},\ft{i}$ act 
on $B\ot B'$ by

\begin{eqnarray*}
\et{i}(b\ot b')&=&\left\{
\begin{array}{ll}
\et{i} b\ot b'&\mbox{ if }\vphi_i(b)\ge\veps_i(b')\\
b\ot \et{i} b'&\mbox{ if }\vphi_i(b) < \veps_i(b'),
\end{array}\right. \\
\ft{i}(b\ot b')&=&\left\{
\begin{array}{ll}
\ft{i} b\ot b'&\mbox{ if }\vphi_i(b) > \veps_i(b')\\
b\ot \ft{i} b'&\mbox{ if }\vphi_i(b)\le\veps_i(b').
\end{array}\right. 
\end{eqnarray*}
Here $0\ot b'$ and $b\ot 0$ should be understood as $0$. For crystals we
are considering, there exists a unique isomorphism $B\ot B'
\stackrel{\sim}{\rightarrow}B'\ot B$, {\em i.e.} a unique map
which commutes with the action of Kashiwara operators. In particular,
it preserves the weight.

For a crystal $B$ we define its affinization 
$\Aff(B)=\{z^d b\mid d\in\Z,b\in B\}$.
$z$ is called the {\em spectral parameter}. $\Aff(B)$ also admits the crystal
structure by $\et{i}\cdot z^d b=z^{d+\delta_{i0}}(\et{i}b),
\ft{i}\cdot z^d b=z^{d-\delta_{i0}}(\ft{i}b)$. The crystal isomorphism 
$B\ot B'\stackrel{\sim}{\rightarrow}B'\ot B$ is lifted up to a map 
$\Aff(B)\ot\Aff(B')\stackrel{\sim}{\rightarrow}\Aff(B')\ot\Aff(B)$
called the {\em combinatorial $R$-matrix} (combinatorial $R$ for short).
It has the following form:
\begin{eqnarray*}
R\;:\;\Aff(B)\ot\Aff(B')&\longrightarrow&\Aff(B')\ot\Aff(B)\\
z^d b\ot z^{d'} b'&\longmapsto&z^{d'+H(b\ot b')}\bt'\ot z^{d-H(b\ot b')}\bt,
\end{eqnarray*}
where $b\ot b'\mapsto\bt'\ot\bt$ under the isomorphism $B\ot B'
\stackrel{\sim}{\rightarrow}B'\ot B$. $H(b\ot b')$ is called the 
{\em energy function} and determined up to a global additive constant by
\[
H(\et{i}(b\ot b'))=\left\{%
\begin{array}{ll}
H(b\ot b')+1&\mbox{ if }i=0,\ \vphi_0(b)\geq\veps_0(b'),\ 
\vphi_0(\bt')\geq\veps_0(\bt),\\
H(b\ot b')-1&\mbox{ if }i=0,\ \vphi_0(b)<\veps_0(b'),\ 
\vphi_0(\bt')<\veps_0(\bt),\\
H(b\ot b')&\mbox{ otherwise }.
\end{array}\right.
\]

Now we can state

\begin{proposition}[Yang-Baxter equation] \label{prop:YBeq}
The following equation holds on 
$\Aff(B_l)\ot\Aff(B_{l'})\ot\Aff(B_{l''})$:
\[
(R\ot1)(1\ot R)(R\ot1)=(1\ot R)(R\ot1)(1\ot R).
\]
\end{proposition}

Suppose $l \ge k$.
In the main text we fix the values of the energy function $H$
on $B_l \ot B_k$ so that 
\begin{equation}\label{eq:Hnormalize}
H(u_l \ot u_k) = 2\varsigma k\quad \text{ for any } \geh,
\end{equation}
regardless of the rank of $\geh$.
See (\ref{eq:varsigma}) for $\varsigma$.
The values given in {\sc Propositions} 
\ref{prop:hweI} and \ref{prop:hweII} are so normalized.
In fact $\max \{ H(b \ot b') \mid b \in B_l, b' \in B_k \} = 2\varsigma k$ holds, 
and as for the minimum we have 
\begin{align*}
&\min \{ H(b \ot b') \mid b \in B_l, b' \in B_k \}
= k,0,2k,0,k,0,2k \nonumber \\
&\text{ for }
\geh = A^{(1)}_n, A^{(2)}_{2n-1}, A^{(2)}_{2n},
B^{(1)}_n, C^{(1)}_n, D^{(1)}_n, D^{(2)}_{n+1}. 
\end{align*}

\section{\mathversion{bold}Crystal $B_l$}\label{app:Bl}

Let us present the parameterization of the crystal
$B_l$  as a set for each $\geh$.
%

%
\noindent
$\geh=A^{(1)}_{n}:$
\begin{equation*} 
B_l = \{ (x_{1},\dots ,x_{n+1}) \in
\Z^{n+1}  | x_{i} \geq 0,\,
\sum_{i=1}^{n+1}x_{i} =l \}.
\end{equation*}
For $B_1$ we use a simpler notation
\begin{equation*}
(x_{1},\dots ,x_{n+1}) =
\hako{i}  \text{ if }x_i=1,\,\text{others}=0.
\end{equation*}

\noindent
$\geh=A^{(2)}_{2n-1}:$
\begin{equation*} \label{eq:BofA2o}
B_l = \{ (x_{1},\dots ,x_{n},\ol{x}_{n},\dots ,\ol{x}_{1}) \in
\Z^{2n}  | x_{i}, \ol{x}_{i} \geq 0,\,
\sum_{i=1}^{n}(x_{i}+\ol{x}_{i}) =l \}.
\end{equation*}
For $B_1$ we use a simpler notation
\begin{equation}\label{eq:0-nashi}
(x_{1},\dots ,x_{n},\ol{x}_{n},\dots ,\ol{x}_{1}) =
\begin{cases}
\hako{i} & \text{if }x_i=1,\,\text{others}=0, \\
\hako{\ol{i}} & \text{if }\ol{x}_i=1,\,\text{others}=0.
\end{cases}
\end{equation}

\noindent
$\geh=A^{(2)}_{2n}:$
\begin{equation*} \label{eq:BofA2e}
B_l = \{  (x_1,\dots,x_n,\ol{x}_n,\dots,\ol{x}_1) \in
\Z^{2n} \;| \; x_i,\ol{x}_i \ge 0,\, \sum_{i=1}^n (x_i + \ol{x}_i)
\le l \}.
\end{equation*}
For $B_1$ we use a simpler notation
\begin{equation*}
(x_1,\dots,x_n,\ol{x}_n,\dots,\ol{x}_1) =
\begin{cases}
\hako{i} & \text{if }x_i=1,\,\text{others}=0, \\
\hako{\ol{i}} & \text{if }\ol{x}_i=1,\,\text{others}=0, \\
\phi & \text{if }x_i=0,\,\ol{x}_i=0 \text{ for all }i.
\end{cases}
\end{equation*}

\noindent
$\geh=B^{(1)}_{n}:$\\
\begin{equation*} \label{eq:BofB}
B_l = \left\{ (x_{1},\dots ,x_{n},x_{0},\ol{x}_{n},\dots ,\ol{x}_{1}) \in
\Z^{2n} \times\{0,1\}
\left|
\begin{array}{c} x_{0} = 0\; \mbox{or} \; 1, x_{i}, \ol{x}_{i} \geq 0, \\
x_0 + \sum_{i=1}^{n}(x_{i}+\ol{x}_{i}) =l
\end{array}\right.
\right\}.
\end{equation*}
For $B_1$
we use a simpler notation
\begin{equation*}
(x_{1},\dots ,x_{n},x_{0},\ol{x}_{n},\dots ,\ol{x}_{1}) =
\begin{cases}
\hako{i} & \text{if }x_i=1,\,\text{others}=0, \\
\hako{\ol{i}} & \text{if }\ol{x}_i=1,\,\text{others}=0.
\end{cases}
\end{equation*}

\noindent
$\geh=C^{(1)}_{n}:$\\
\begin{equation*} \label{eq:BofC}
B_l = \left\{  (x_1,\dots,x_n,\ol{x}_n,\dots,\ol{x}_1) \in
\Z^{2n} \Biggm|
\begin{array}{l}
x_i,\ol{x}_i \ge 0, \\
l \ge \sum_{i=1}^n
(x_i + \ol{x}_i) \in l - 2 \Z
\end{array}
\right\}.
\end{equation*}
For $B_1$ we
use the notation (\ref{eq:0-nashi}).

\noindent
$\geh=D^{(1)}_{n}:$\\
\begin{equation*} \label{eq:BofD}
B_l = \left\{ (x_{1},\dots ,x_{n},\ol{x}_{n},\dots ,\ol{x}_{1}) \in
\Z^{2n}  \left|
\begin{array}{c}
x_{n} = 0\; \mbox{or} \;\ol{x}_{n} =0,\; x_{i}, \ol{x}_{i} \geq 0,\\
\sum_{i=1}^{n}(x_{i}+\ol{x}_{i}) =l
\end{array}\right.
\right\}.
\end{equation*}
For $B_1$
we use the notation (\ref{eq:0-nashi}).

\noindent
$\geh=D^{(2)}_{n+1}:$\\
\begin{equation*} \label{eq:BofD2}
B_l = \left\{  (x_1,\dots,x_n,x_0,\ol{x}_n,\dots,\ol{x}_1) \in
\Z^{2n} \times \{0,1\} \Biggm|
\begin{array}{l}
x_0=\mbox{$0$ or $1$},x_i,\ol{x}_i \ge 0, \\
x_0 + \sum_{i=1}^n (x_i + \ol{x}_i) \le l
\end{array}
\right\}.
\end{equation*}
For $B_1$
we use a simpler notation
\[
B_1 \ni (x_1,\dots,x_n,x_0,\ol{x}_n,\dots,\ol{x}_1) =
\begin{cases}
\hako{i} & \text{if }x_i=1,\,\text{others}=0, \\
\hako{\ol{i}} & \text{if }\ol{x}_i=1,\,\text{others}=0, \\
\phi & \text{if }x_i=0,\,\ol{x}_i=0 \text{ for all }i.
\end{cases}
\]

\section{\mathversion{bold} Crystal structure of $B_\natural$}\label{section:B-natural}

We review the crystal structure of $B_\natural$ for Type I and II, {\it i.e.}, 
for $\geh=A^{(2)}_{2n-1},B^{(1)}_n,D^{(1)}_n,C^{(1)}_n$. Define an index set $J$ and
a linear order $\prec$ on $J$ by
\begin{eqnarray*}
J&=&\{1\prec2\prec\cdots\prec n\prec0\prec\ol{n}\prec\cdots\prec\ol{2}\prec\ol{1}\}
\quad\mbox{ for }\geh=B^{(1)}_n,\\
&=&\{1\prec2\prec\cdots\prec n\prec\ol{n}\prec\cdots\prec\ol{2}\prec\ol{1}\}
\quad\mbox{ for }\geh=A^{(2)}_{2n-1},D^{(1)}_n,C^{(1)}_n.
\end{eqnarray*}
For $\geh$ of Type I or II, let $\gehb$ be the classical part  
specified in the beginning of Section \ref{subsection:plan}.
As a set $B_\natural$ is given by 
\begin{eqnarray*}
B_\natural&=&B(\Lab_2)\oplus B(0)\quad\mbox{for Type I }
(\geh=A^{(2)}_{2n-1},B^{(1)}_n,D^{(1)}_n),\\
&=&B(\Lab_2)\quad\mbox{for Type II }(\geh=C^{(1)}_n),\\
B(0)&=&\{\phi\},\\
B(\Lab_2)&=&\left\{{\alpha\choose\beta}\left|
\begin{array}{l}
\alpha,\beta\in J\mbox{ satisfying conditions}\\
\mbox{(1) and (2) below}
\end{array}\right.\right\}.
\end{eqnarray*}
\begin{itemize}
\item[(1)] $\alpha\prec\beta$ \quad if $\ol{\geh}=C_n$,\\
$\alpha\prec\beta$ or $(\alpha,\beta)=(0,0)$ \quad if $\ol{\geh}=B_n$,\\
$\alpha\prec\beta$ or $(\alpha,\beta)=(\ol{n},n)$ \quad if $\ol{\geh}=D_n$.
\item[(2)] $(\alpha,\beta)\neq(1,\ol{1})$.
\end{itemize}
The labeling of the elements of $B(\Lab_2)$ is taken {}from \cite{KN}.
Needless to say, the actions of the crystal operators $\et{i},\ft{i}
(i=1,2,\cd,n)$ on $B(\Lab_2)$ agree with those in \cite{KN}. 
On $B(0)$ they are defined by
\[
\et{i}\phi=\ft{i}\phi=0\quad(i=1,2,\cd,n).
\]
Therefore we are left to show the actions of $\et{0},\ft{0}$, which are given
as follows:
\par\noindent
For Type I,
\begin{eqnarray*}
&&{\ol{2}\choose\ol{1}}\stackrel{0}{\longrightarrow}\phi
\stackrel{0}{\longrightarrow}{1\choose2},\\
&&{\alpha\choose\ol{2}}\stackrel{0}{\longrightarrow}{1\choose\alpha}
\quad\mbox{if }\alpha\neq1,2,\\
&&{\alpha\choose\ol{1}}\stackrel{0}{\longrightarrow}{2\choose\alpha}
\quad\mbox{if }\alpha\neq2,\ol{2},\\
&&{\alpha\choose\beta}\stackrel{0}{\longrightarrow}0\quad\mbox{in the other cases}.
\end{eqnarray*}
\par\noindent
For Type II,
\begin{eqnarray*}
&&{\alpha\choose\ol{1}}\stackrel{0}{\longrightarrow}{1\choose\alpha}
\quad\mbox{if }\alpha\neq1,\\
&&{\alpha\choose\beta}\stackrel{0}{\longrightarrow}0\quad\mbox{in the other cases}.
\end{eqnarray*}
As usual, $b\stackrel{i}{\longrightarrow}b'$ means $\ft{i}b=b'\Longleftrightarrow
b=\et{i}b'$.


\end{document}